\documentclass[journal]{new-aiaa} 
\usepackage[utf8]{inputenc}

\usepackage{graphicx}
\usepackage[justification=centering]{caption}
\usepackage{amsmath,amsfonts}
\usepackage[version=4]{mhchem}
\usepackage{siunitx}
\usepackage{longtable,tabularx,changepage}
\usepackage{color,soul}
\usepackage{subcaption,float}
\usepackage{hyperref}
\usepackage{enumitem}
\usepackage[flushleft]{threeparttable}
\usepackage{rotating}
\usepackage[section]{placeins}
\usepackage{flafter}
\usepackage{lineno}
\usepackage{cancel,bm}
\usepackage{soul}
\usepackage[export]{adjustbox}
\allowdisplaybreaks

\newcounter{subeq}
\renewcommand{\thesubeq}{\theequation\alph{subeq}}
\newcommand{\newsubeqblock}{\setcounter{subeq}{0}\refstepcounter{equation}}
\newcommand{\mysubeq}{\refstepcounter{subeq}\tag{\thesubeq}}
\newcommand{\obj}{\mathcal{J}}
\newcommand{\hquad}{\hspace{0.5em}} 

\def\code#1{\texttt{#1}}

\title{Event-Driven Network Model for Space Mission Optimization with High-Thrust and Low-Thrust Spacecraft \footnote{This paper is a substantially revised version of the paper AAS 18-337 presented at the AAS/AIAA Astrodynamics Specialist Conference, Snowbird, UT, August 19-23, 2018 \cite{jag_c3}.}}

\author{Bindu B. Jagannatha \footnote{PhD, Department of Aerospace Engineering, Urbana, IL, 61801.}}
\affil{University of Illinois at Urbana-Champaign, Urbana, IL, 61801}
\author{Koki Ho \footnote{Assistant Professor, Daniel Guggenheim School of Aerospace Engineering, Atlanta, GA 30332, AIAA Member.}}
\affil{Georgia Institute of Technology, Atlanta, GA, 30332}

\begin{document}
\maketitle

\small
\begin{abstract}
Numerous high-thrust and low-thrust space propulsion technologies have been developed in the recent years with the goal of expanding space exploration capabilities; however, designing and optimizing a multi-mission campaign with both high-thrust and low-thrust propulsion options are challenging due to the coupling between logistics mission design and trajectory evaluation. Specifically, this computational burden arises because the deliverable mass fraction (i.e., final-to-initial mass ratio) and time of flight for low-thrust trajectories can can vary with the payload mass; thus, these trajectory metrics cannot be evaluated separately from the campaign-level mission design. To tackle this challenge, this paper develops a novel event-driven space logistics network optimization approach using mixed-integer linear programming for space campaign design. An example case of optimally designing a cislunar propellant supply chain to support multiple lunar surface access missions is used to demonstrate this new space logistics framework. The results are compared with an existing stochastic combinatorial formulation developed for incorporating low-thrust propulsion into space logistics design; our new approach provides superior results in terms of cost as well as utilization of the vehicle fleet. The event-driven space logistics network optimization method developed in this paper can trade off cost, time, and technology in an automated manner to optimally design space mission campaigns.
\end{abstract}

\section*{Nomenclature}

{\renewcommand\arraystretch{1.0}
	\noindent\begin{longtable*}{@{}l @{\quad=\quad} l@{}}		
		$\mathbb{A}_\text{H}$ & set of holdover arcs\\
		$\mathbb{A}_\text{T}$ & set of transportation arcs\\
		$\bm{B}$  & mass transformation matrix\\
		$\bm{C}$  & concurrency matrix\\
		$\bm{c}$  & cost coefficient vector\\
		$\bm{d}$  & demand/supply vector\\
		$\mathbb{E}$ & set of events\\
		$e$ & event index\\
		$\bm{F}$ & parameter matrix  for time-related constraints	\\
		$i$  & node index\\
		$j$  & node index\\
		$\bm{M}$  & mass matrix\\		
		$\mathbb{N}$ & set of nodes \\
		$p$  & parameter for affine/linear approximation \\
		$q$  & parameter for affine/linear approximation \\
		$\bm{s}$ & parameter vector for concurrency constraints\\		
		$\bm{t}$ & parameter vector for time-related constraints\\		
		$v$ & vehicle index\\
		$\bm{x}$ & commodity flow vector\\
		$y$ & total commodity flow\\
		$\lambda,\lambda'$ & auxiliary variable vector for piecewise linear approximation \\ 
		$\bm{\tau}$ & auxiliary variable vector for time-related constraints	
\end{longtable*}}

\section{Introduction}
The complexity of space mission design is increasing as newer technologies such as solar electric propulsion (SEP), in-situ resource utilization (ISRU), in-space manufacturing, on-orbit propellant depots, etc. mature into space-ready states. These emerging technologies provide multiple avenues for long-term human space exploration through multi-mission campaigns aimed at incrementally building up capabilities. In this pursuit of sustaining space exploration goals, space logistics tools can help campaign designers make architectural decisions by balancing all associated costs of using different technologies through rigorous mathematical modeling. Such tools can conduct top-level trades in an automated manner and evaluate potential approaches to meeting long-term strategic objectives. More importantly, by exploring a broad range of logistics strategies, solutions may emerge that expertise-based designs could miss. 

Space logistics optimization, thus, is a growing field of study for reducing the costs of long-duration campaigns involving deployment and utilization of space infrastructure. Trajectory analysis plays a key role in deciding the costs of missions that comprise the campaign. Traditional trajectory design methods are often confined to analyzing individual missions, whether using high-thrust or low-thrust propulsion options, and do not typically consider the architectural aspects of the in-space network. On the other hand, the determination of network architecture and mission sequence often relies on the expertise of mission designers (e.g., OASIS \cite{oasis2002}). More recently, space logistics planning tools have been developed to guide such decisions through mathematical modeling. For example, optimization of multi-commodity network flow using combinatorial heuristic algorithms was suggested by Taylor et al. \cite{taylor1,taylor2}. This method was improved by Ishimatsu et al. through the development of the generalized multi-commodity network flow (GMCNF) model \cite{static1,static2}. These space logistics methods have been extended to campaign-level mission design framework using a time-expanded network by Ho et al \cite{dynamic1,dynamic2}. These frameworks approach the problem of campaign design through a logistics-driven perspective. Other works such as the EXAMINE framework \cite{examine} and the graph-theory-based space system architectures model developed by Arney et al. \cite{arney} can compare and/or optimize user-input system architecture scenarios. 

A critical limitation of the conventional space logistics design methods referenced above is that they are unable to account for the use of low-thrust vehicles for transportation.  
More specifically, the conventional methods assumed a decoupling between logistics mission design and trajectory evaluation, where the mission design model takes pre-computed trajectory metrics such as the final-to-initial mass ratio (or $\Delta$v) and time of flight for each arc as inputs and then optimize the mission architecture including the logistics flows of propellant, payload, and other commodities. This approach is effective for a campaign with a fleet of high-thrust spacecraft, and responded well to the background that many past space exploration campaigns considered only high-thrust spacecraft options (e.g., chemical propulsion). 
However, given the significant advancements made in low-thrust propulsion technology (e.g., SEP) in recent years, there is a growing interest in the question of how low-thrust spacecraft can be optimally integrated into a space exploration campaign. Thus, we need a mathematical framework that can optimally identify when/where to use high-thrust vs. low-thrust spacecraft options in a campaign.
This research question is challenging because of the coupling between logistics mission design and trajectory evaluation for low-thrust spacecraft. Namely, the final-to-initial mass ratio and time of flight for low-thrust trajectories can vary with the payload mass, and thus cannot be pre-computed for mission design. We need the evaluation of these trajectory metrics to be integrated into the logistics mission design process so that it is coupled with the decisions for the flows of propellant, payload, and other commodities. As reviewed in detail later in Section~\ref{sec:background}, the conventional space logistics methods cannot effectively handle the complexity due to this coupling.

In response to this background, the aim of this paper is to develop a unified space logistics method to design space exploration campaigns with both high-thrust and low-thrust trajectory technologies. We expect that this contribution will broaden the tradespace of astrodynamics studies and enable the space mission designers to explore architectural options that use different propulsion technologies cooperatively.
The new formulation developed in this work is demonstrated in a case study where cislunar infrastructure is deployed for in-space propellant resupply of crew missions to the lunar surface. While the crew missions themselves use high-thrust vehicles, the resupply propellant tanks can be predeployed to cislunar locations using a fleet of cargo delivery tugs consisting of high-thrust and low-thrust variants. These tugs provide cost savings by exploiting cheaper (albeit slower) routes; this also decouples the crew flight time from the total campaign duration. The state-of-the-art method to analyze this problem is to use a general optimization approach with genetic algorithms \cite{jag_j1}. However, as demonstrated later in this paper, the formulations proposed in the existing literature require a number of simplifying assumptions and do not explore the tradespace efficiently. Instead, this paper develops a novel network modeling approach using mixed-integer linear programming to tackle this problem. The proposed event-driven space logistics network model concurrently optimizes three aspects of the problem -- the best route for the crew missions within the given network if refuel tanks can be made available, the most strategic locations for these refuel tanks to be positioned, and the best combination of available high-thrust and low-thrust cargo tugs to predeploy these refuel tanks. The developed formulation incorporates surrogate models for trajectory evaluation and uses them to concurrently optimize the campaign logistics and trajectory selection. The surrogate models for our case study are either obtained from existing trajectory models for high-thrust propulsion or derived from applying the computationally efficient Q-law low-energy low-thrust trajectory design method \cite{jag_c2,jag_j2} to the current problem.

Mathematical modeling and optimization of space logistics developed in this paper can support campaign-level architectural decisions where different propulsion technologies can be traded off to achieve overall exploration goals. The decisions made in this manner have a higher potential of sustaining long-term human and robotic space exploration than traditional mission-wise analyses.

\section{Background and Motivation}
\label{sec:background}
This work builds on the underlying mixed-integer linear programming (MILP) formulation of the space logistics model by introducing a new concept of event-driven networks to incorporate both high-thrust and low-thrust trajectories. This section briefly reviews the conventional space logistics model, and presents an overview of the challenges faced by its existing derivatives when dealing with low-thrust trajectories.

\subsection{Review of the Conventional Space Logistics Optimization Methods}
The conventional space logistics model is based on the GMCNF formulation \cite{static1,static2}. It begins by modeling the space exploration map as a network, with nodes corresponding to the candidates for surface destinations, orbits, and staging locations. Arcs connect these nodes, allowing materials to flow between them. Thus, each of these arcs has parameters associated with it, like the payload/vehicle combination, time, technology used on the arc, etc. In the GMCNF framework, these parameters are modeled as commodity flows, where commodities include crew, vehicles, and propellant, along with required structural masses of the propulsive elements. The cost of transporting commodities across the network to meet the demand at all nodes is minimized while conforming to network flow principles such as mass balance, flow transformation, and flow concurrency. This framework combines the classical concepts of generalized network flow and multi-commodity network flow, and applies them to the design of space exploration campaigns.

\begin{figure}[htb]
\centering
\includegraphics[trim={0 0 11cm 0},clip,width=0.4\textwidth]{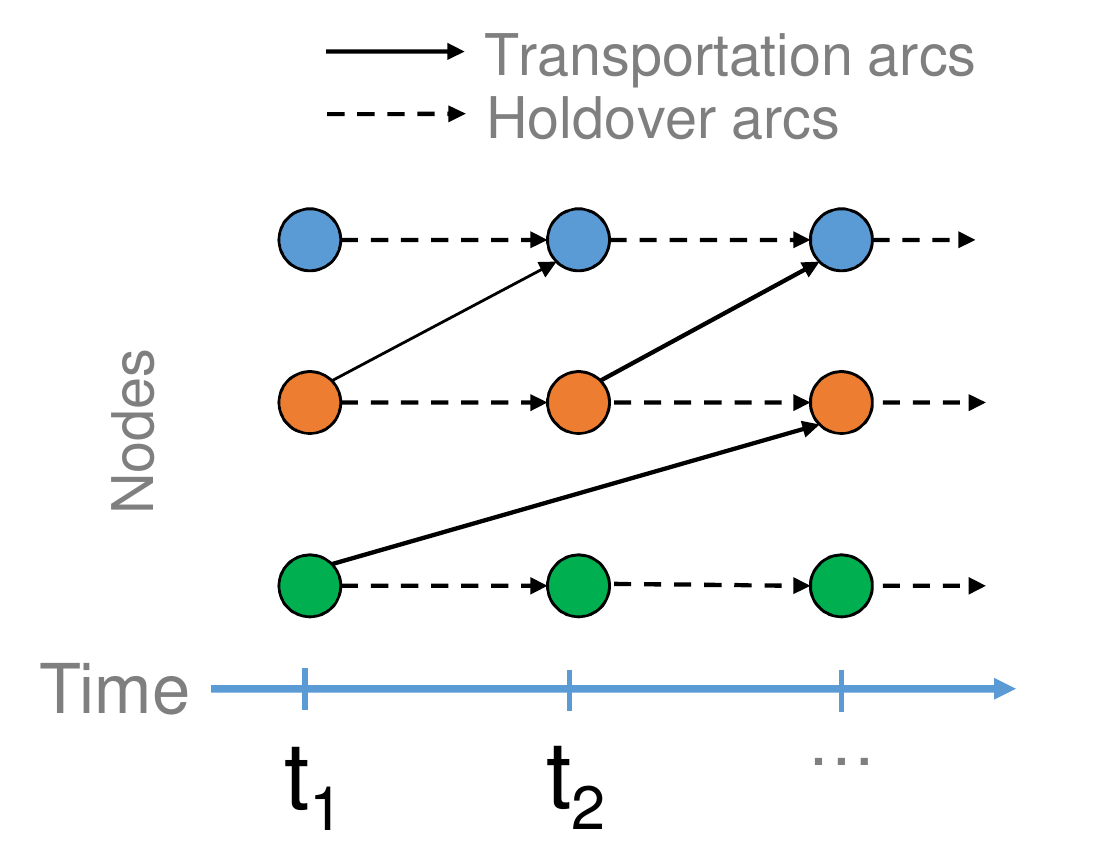}
\caption{Time-expanded network flow model\label{fig:timex}}
\end{figure}

The major drawback of the above original GMCNF model is its inability to deal with temporal behavior correctly, which results in time paradoxes and infeasible flow-generation loops for campaign-level analysis. These shortcomings related to the static nature of the original space logistics model were addressed by its time-expanded variant \cite{dynamic1}. Figure~\ref{fig:timex} depicts this time-expanded network model -- the time dimension is assimilated by duplicating the underlying static network at discrete time steps separated by known time spans. This notional static network uses \textit{transportation arcs} to facilitate mass flow between different nodes, while \textit{holdover arcs} that connect the same node across time steps are used to incorporate the concept of stock. These arcs only represent the possibility of flow, and need not be all utilized in the actual solution. 

This updated time-expanded network model is a powerful formulation that is capable of lifecycle optimization of mission sequences within a campaign. However, the time-expanded network model and its bi-scale variant \cite{dynamic2} are based on the assumption of fixed arc lengths, i.e., the flight time across each arc is known before setting up the network and is independent of the commodity flow. This assumption is incompatible with certain space transportation options, like the low-thrust trajectories, as discussed next.

\subsection{Challenge of Space Trajectory Evaluation}
As reviewed previously, space transportation options can broadly be classified into high-thrust and low-thrust categories, and the key performance metrics associated with using these modes of space transportation are two-fold: final-to-initial mass ratio and time of flight. These metrics for a trajectory with a given high-thrust vehicle can be treated as constant (i.e., independent of the payload mass) because the required in-space propulsive maneuvers are modeled as impulsive burns. In contrast, low-thrust transfers require the engine to be continuously powered over extended periods of time due to the low levels of thrust acceleration provided by such propulsion systems. Thus, the performance metrics for a low-thrust trajectory are coupled with the thrust-to-mass ratio of the propulsive vehicle; this indicates that, in general, for a given low-thrust vehicle specification and trajectory option (e.g., constant thrust, constant $I_\text{sp}$), the final-to-initial mass ratio and/or time of flight can depend on the payload mass to be delivered. Particularly, from the space logistics design perspective, this means that the transit time for a space transportation arc (or the arc length) when using low-thrust propulsion can be dependent on the commodity flow on that arc. 

Given this background, the reader can notice that the use of discrete time steps in the time-expanded network model does not align well with the nonconstant nature of low-thrust transportation arc lengths. A similar issue of flow-dependent arc lengths can also be observed in terrestrial logistics applications such as road traffic and communication networks, where the amount of time required to transit an arc increases as more traffic enters it. Conventional logistics frameworks have treated this problem through the use of multiple arcs of different capacities connecting the same pair of nodes in a time-expanded network, with each arc representing a given discrete transit-time option \cite{fdal1,fdal3}. 
However, for the existing approaches to model the vehicle flow with an accurate arrival time in a general network optimization problem, we need a fine time-step discretization for the time-expanded network. For the space logistics case we are considering, the flight time of a low-thrust trajectory can take a wide range of values (from days to years) depending on the payload mass, and a fine time-step discretization (e.g., in days) can result in an intractably large time-expanded network and thus cannot be solved practically.

Furthermore, low-thrust trajectory design is a challenging problem in itself, with an entire branch of astrodynamics research dedicated to it. This is because evaluating the performance of any low-thrust transfer requires determining the continuous time history of the thrust control vector for the entire trajectory. Analytical solutions to low-thrust transfers exist only in a small class of academic problems, and thus computationally-expensive nonlinear optimization methods are often employed to calculate such trajectories. Even after deriving performance models for low-thrust transportation through complex trajectory analyses, these models cannot yet be incorporated into space logistics frameworks at a campaign-level in an efficient manner. By dealing with this challenge, the field of trajectory design can be more fully integrated into space logistics planning -- this is the context that motivates the ideas developed in this paper.

\section{Methodology}
The framework introduced here for incorporating flow-dependent transit times into space logistics modeling frameworks is based on expanding the static network model across the time dimension by using \textit{event-based steps}.  Each event ``layer'' within the network spans some time period that is assigned internally to satisfy the time constraint on the total campaign duration. This event-driven expansion of the static network accommodates the trajectory performance models in a way that tradeoffs between available propulsion technologies can be conducted internally. 

The current section presents the formulation of the event-driven network flow model, which is the newly developed alternative to the existing time-expanded network model for optimizing space logistics problems involving flow-dependent arc transit times. First, the mixed-integer linear programming (MILP) formulation used here to arrive at the optimal logistics solution is detailed, and then the techniques to treat nonlinearities arising from trajectory design surrogate models are addressed. 

\subsection{Event-Driven Network Flow Model}
\begin{figure}[ht!]
\centering
\includegraphics[trim={11cm 0 0 0},clip,width=0.4\textwidth]{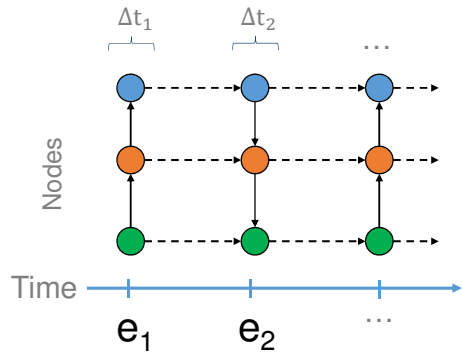}
\caption{Event-driven network flow model\label{fig:eventex}}
\end{figure}

The new space logistics framework approaches the issue of flow-dependent transit times by duplicating the static network at ``events'', instead of at time steps. Each copy of the underlying static network, or ``layer'', encompasses a pre-defined event and allows commodity flow between nodes via intra-layer transportation arcs. Unlike in the time-expanded network model, these transportation arcs do not connect nodes across layers; instead, (inter-layer) holdover arcs are used to connect copies of the same node across all the events. Figure~\ref{fig:eventex} helps contrast this event-driven expansion of the static network with the previous time expansion. 

While each layer is a copy of the underlying static network, some arcs are marked active or inactive depending on the event assigned to the layer. Active arcs simply denote the possibility of commodity flow and it may be that not all of them are utilized in the final solution. Inactive arcs should be removed from the event layer as no commodity flow is allowed across them. The distribution of events is campaign-specific, but must follow these basic rules: 
\begin{itemize}
    \item The active arcs in an event layer should not form a cyclic graph. This unidirectionality within a single event layer prevents undesirable flow-generation loops. 
    \item Each event layer can consist of one or more transportation arcs as long as no commodity exchange is allowed among the flowing vehicles within that event. This is because the exact time of arrival of commodities at nodes is not tracked, and thus commodity exchange cannot be controlled. Note that this also indicates that multiple routes may connect the same pair of nodes within the same layer as long as they do not interact with each other. 
    \item The time length assigned to each event should correspond to the greatest of the longest utilized directed paths between any (feasible) pair of nodes in the event layer. A utilized path is defined as a continuous sequence of connected active arcs with nonzero commodity flow.
\end{itemize}

\noindent Thus, a notional static network can be devised in this way to represent the dynamic behavior of the campaign. Using the formulation derived in the next section, the time span of each event can be adjusted internally as required to meet the desired time constraints. The use of events is explained further within the case study, after the mathematical formulation for the event-driven network model is derived. 

\subsection{Mathematical Formulation}
\begin{figure}[ht!]
\captionsetup[subfigure]{justification=centering}
\centering
\begin{subfigure}[t]{0.32\textwidth}
\centering
\includegraphics[width=\textwidth]{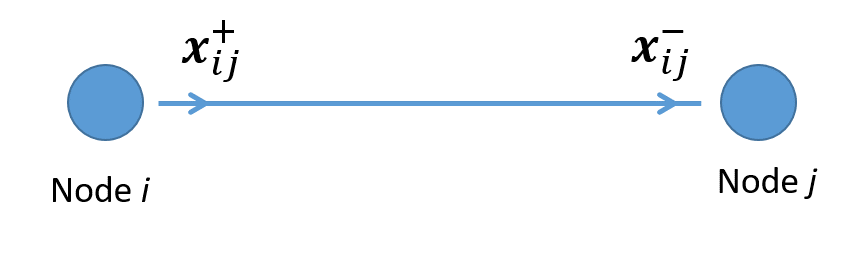}
\caption{Single transportation arc\label{fig:arc1}}
\end{subfigure}
~ 
\begin{subfigure}[t]{0.32\textwidth}
\centering
\includegraphics[width=\textwidth]{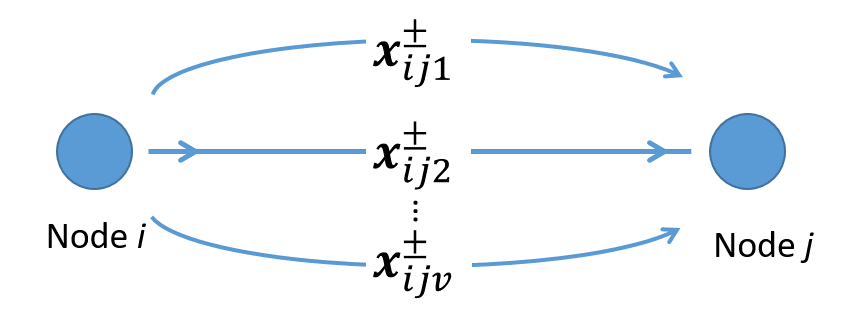}
\caption{Multi-graph transportation arc\label{fig:arc2}}
\end{subfigure}
~
\begin{subfigure}[t]{0.32\textwidth}
\centering
\includegraphics[width=\textwidth]{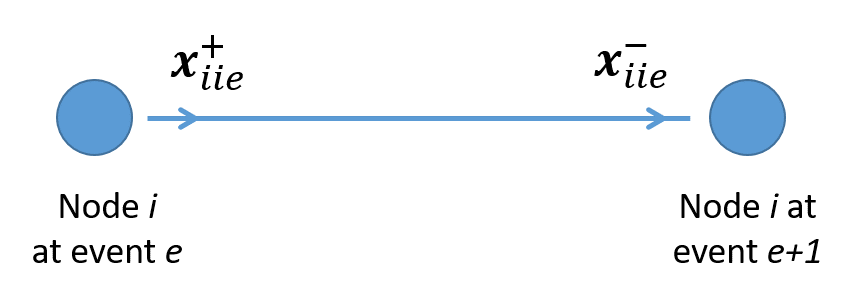}
\caption{Holdover arc\label{fig:arc3}}
\end{subfigure}
\caption{Arcs used in network graph formulation.\label{fig:arc}}
\end{figure}

Consider a static graph consisting of a set of nodes $\mathbb{N}_\text{static}$, connected by a set of directed arcs $\mathbb{A}_\text{static}$, which allow transportation of multiple commodities across it. For each arc from the node $i$ to node $j$, the multi-commodity flow can be split into outflow $\bm{x}^+_{ij}$ and inflow $\bm{x}^-_{ij}$ as shown in Figure~\ref{fig:arc1}. Additionally, a cost coefficient $\bm{c}^+_{ij}$ is assigned to the outflow from each node. Within this static graph, different technology options along each arc can be represented by the multi-graph (Figure \ref{fig:arc2}), where multiple arcs connect the same pair of end nodes. This concept of multi-graph is used to denote the choice of discrete alternatives for propulsion in the current work. The arc index thus assumes the form ($i,j,v$), where $v\in \mathbb{V}$ with $\mathbb{V}$ denoting the set of all propulsive vehicles available in the analysis. 

Now, consider the dynamic event-driven network model where expansion in the time dimension is conducted by repeating the static network at event-steps instead of discrete time steps. These events form a pre-determined set $\mathbb{E}$ and are marked by the index \textit{e} where $e\in\mathbb{E}$. Due to this expansion of the static model into the time dimension, each node has an index ($i,e$) meaning the node at event-step $e$. Each transportation arc now has an index ($i,j,v,e$) and is used to represent a flow of commodities from node $i$ to node $j$ by using propulsion provided by vehicle $v$ within the event-step $e$. Holdover arcs, which facilitate the flow of commodities between these static layers, connect the copies of the same node $i$ across event steps $e$ and $e+1$ and hence are denoted by the index $(i,i,e)$. As shown in Figure \ref{fig:arc3}, holdover arcs are always directed forward in time, which means that a holdover arc originating in event $e$ cannot connect to its node copy in event $e-1$. 
The set of nodes in this event-driven network is denoted as $\mathbb{N}$; the set of transportation arcs and that of holdover arcs are denoted as $\mathbb{A}_\text{T}$ and $\mathbb{A}_\text{H}$, respectively. 

We are interested in optimizing the multi-commodity flow $\bm{x}^\pm_{ijve}$ (for transportation arcs) and $\bm{x}^\pm_{iie}$ (for holdover arcs) over this event-driven network.
Due to the multi-commodity nature, these variables are vectors; if $k$ commodities flow across the network graph, then $\bm{x}^\pm_{ijve}$ and $\bm{x}^\pm_{iie}$ are $k$-by-1 vectors. These commodities can be integer-valued or continuous in nature; for example, crew and vehicles should be an integer-valued commodity. A diagonal matrix $\bm{M}$ can be used to convert these integer-valued commodity variables into corresponding masses as required.
We also define $\bm{x}^+$ 
as the consolidated vector of $\bm{x}^+_{ijve}$ 
across the entire set of the transportation arcs $(i,j,v,e)$. 
Finally, we define an auxiliary variable vector $\bm{\tau}$, which is needed for the time-related constraints as discussed later.

Using the above developed notation, the event-driven network model can be expressed by the following formulation:
\begin{equation}
\mathrm{Minimize\quad } \obj =  \sum_{e\in \mathbb{E}} \sum_{(i,j,v,e)\in \mathbb{A}_\text{T}} {\bm{c}^+_{ijv}}^T \bm{x}^+_{ijve} + \sum_{e\in \mathbb{E}} \sum_{(i,i,e)\in \mathbb{A}_\text{H}} {\bm{c}^+_{ii}}^T \bm{x}^+_{iie}\hquad, \label{eq:cost}
\end{equation}
\noindent subject to the following constraints:
\begin{align}
\sum_{\substack{(j,v):\\(i,j,v,e)\in\mathbb{A}_\text{T}}} \bm{x}^+_{ijve} - \sum_{\substack{(j,v):\\(j,i,v,e)\in\mathbb{A}_\text{T}}} \bm{x}^-_{jive} + \bm{x}^+_{iie} &- \bm{x}^-_{ii(e-1)} \leq \bm{d}_{ie} \quad \forall\ (i,e) \in \mathbb{N}\label{eq:balance}\hquad,\\
\newsubeqblock
\mysubeq \bm{B}_{ijv}\bm{x}^+_{ijve} = \bm{x}^-_{ijve} &\qquad \forall\ (i,j,v,e) \in \mathbb{A}_\text{T} \label{eq:xformation1} \hquad,\\
\mysubeq \bm{B}_{ii}\bm{x}^+_{iie} = \bm{x}^-_{iie} &\qquad \forall\ (i,i,e) \in \mathbb{A}_\text{H} \label{eq:xformation2}\hquad,\\
\newsubeqblock
\mysubeq \bm{C}^+_{ijv}\bm{x}^+_{ijve} \leq \bm{s}^+_{ijv} &\qquad \forall\ (i,j,v,e) \in \mathbb{A}_\text{T}  \label{eq:concur1}\hquad,\\
\mysubeq \bm{C}^+_{ii}\bm{x}^+_{iie} \leq \bm{s}^+_{ii} &\qquad \forall\ (i,i,e) \in \mathbb{A}_\text{H}  \label{eq:concur2}\hquad,\\
\newsubeqblock 
\mysubeq \bm{x}^\pm_{ijve} \geq \bm{0}_{k\times 1} &\qquad \forall\ (i,j,v,e) \in \mathbb{A}_\text{T}  \label{eq:bound1} \hquad,\\
\mysubeq \bm{x}^\pm_{iie} \geq \bm{0}_{k\times 1} &\qquad \forall\ (i,i,e) \in \mathbb{A}_\text{H}  \label{eq:bound2} \hquad,\\
\bm{F}\begin{bmatrix}
\bm{x}^+\\
\bm{\tau}
\end{bmatrix} \leq \bm{t} \label{eq:tof_con}\hquad. 
\end{align}

Equation~\eqref{eq:cost} expresses the objective function that should be minimized for efficient campaign design, where $\bm{c}^+_{ijv}$ and $\bm{c}^+_{ii}$ are the cost coefficients. Its first term shows the cost of transporting commodities across the network and its second term represents the cost of holding stock at nodes across layers.

The first constraint, Equation~\eqref{eq:balance}, deals with mass balance by ensuring that the demand at every node is met under the supply conditions. Each node $i$ at event-step $e$ also has an associated demand vector $\bm{d}_{ie}$, which represents the demand or supply of each commodity at that node within that event layer. Demand is indicated by negative components in $\bm{d}_{ie}$, whereas supply is indicated by positive components. 

Equation~\eqref{eq:xformation1} uses the $k$-by-$k$ mass transformation matrix $\bm{B}_{ijv}$ to show gains/losses and account for commodity transformation. Similarly, the commodity gain/loss/transformation over holdover arcs is expressed by Equation~\eqref{eq:xformation2}. If no mass transformation is expected over a certain arc, then the corresponding mass transformation matrix ($\bm{B}_{ijv}$ or $\bm{B}_{ii}$) should be set to equal to an identity matrix $\bm{I}_{k\times k}$. If multiple types of transformations occur over an arc, the corresponding constraint over it can be aggregated by serially multiplying the $\bm{B}$ matrices in an order matching the sequence of transformation processes across that arc. Examples of flow transformations include propellant consumption for transportation, generation of resources through ISRU, consumption of raw materials during in-space additive manufacturing, etc. 

Concurrency constraints are taken care of by Equation~\eqref{eq:concur1} over transportation arcs and by Equation~\eqref{eq:concur2} over holdover arcs. If $n$ such constraints have to be enforced, then the matrix $\bm{C}^+_{ijv}$ (and $\bm{C}^+_{ii}$) would be of the size $n$-by-$k$ and $\bm{s}^+_{ijv}$ (and $\bm{s}^+_{ii}$) a vector of size $n$-by-$1$. These constraints vary across arcs depending on the layer that they belong to, and are listed for the case study in detail later. An example concurrency constraint is the requirement that the amount of propellant being carried across a transportation arc does not exceed the propulsive vehicle's fuel carrying capacity. 

Equations~\eqref{eq:bound1}--\eqref{eq:bound2} represent the nonnegative nature of flows across all arcs. Note that, together with the concurrency constraints in Equations~\eqref{eq:concur1}--\eqref{eq:concur2}, we can force a subset of the variables to be zero; this constraint can be used to prevent the flow of certain commodities over arcs within particular event layers. 

Finally, the last constraint in Equation~\eqref{eq:tof_con} can be expanded to include multiple time-related constraints -- it represents specific time-related bounds that need to be respected by the campaign. 
These constraints are appended to the formulation as a means of exploring the cost-versus-time tradespace for the campaign and represent a new feature that does not exist in the static and the time-expanded space logistics network models. As an example, this type of constraint could be used to constrain the lengths of low-thrust transportation arcs -- the flight time of a low-thrust propulsion spacecraft is dependent on the total mass entering the arc. Another example would be to constrain the total time span of event/s depending on the functional breakdown of events. An auxiliary variable vector $\bm{\tau}$ is introduced to accommodate the time-related restrictions on any given campaign. For example, this vector may contain the time spans of the campaign event layers, so that the time-related constraints may be applied in a linear form. The use of these constraints is demonstrated in the case study campaign in Sections~\ref{sec:case}--\ref{sec:results}.

Although the above constraints represent the complete formulation, for simplicity of notation in the later analysis, we introduce an additional auxiliary variable $y^\pm_{ijve}$ with an additional constraint:
\begin{align}
y^\pm_{ijve} = \bm{1}^T \bm{M}\bm{x}^\pm_{ijve} &\qquad \forall\ (i,j,v,e) \in \mathbb{A}_\text{T} \label{eq:tot_mass}
\end{align}
The $k$-by-1 vector $\bm{1}$ has all its components set to unity and is used here with mass matrix $\bm{M}$ to obtain the sum of the mass of the components of $\bm{x}^\pm_{ijve}$. The resulting variable $y^\pm_{ijve}$ is useful in formulating the constraints that are dependent on this sum total of commodity flow across an arc, instead of the flow of some individual commodity.

Note that all relationships in Equations~\eqref{eq:cost}--\eqref{eq:tot_mass} are linear, and thus the space logistics problem can be effectively optimized by using mixed-integer linear programming (MILP) solvers with a given optimality tolerance. The proposed approach to integrate nonlinear low-thrust trajectory models into the MILP-based space logistics formulation is part of our unique contribution and is introduced in more detail in the next subsection (Section~\ref{subsec:flow-dependent-cost}). The decision variables available to the MILP optimizer are $\left\lbrace \bm{x}^\pm_{ijve}, \hquad y^\pm_{ijve}, 
\hquad \bm{x}^\pm_{iie}, \hquad \bm{\tau} \right\rbrace$.
Additionally, as discussed next, extra variables needed to linearize nonlinear relationships between existing arc parameters are also included in the decision variables.

\subsection{Treating Flow-Dependent Trajectory Performance Metrics \label{subsec:flow-dependent-cost}}
The primary contribution of the event-driven network model is that it provides a way to treat (low-thrust) trajectory performance metrics that are dependent on the arcs' mass flow in the MILP space logistics problem formulation. We propose using surrogate models, derived from the trajectory analysis, to capture the relationship between these trajectory performance metrics (i.e., final-to-initial mass ratio and time of flight) and the arcs' initial mass, and incorporating them into the event-driven space logistics optimization formulation. Two forms of the surrogate model can be handled within the proposed framework -- affine (which includes linear) and nonlinear (through piecewise linear approximation). 

This section introduces how these general surrogate models can be incorporated into the MILP formulation through reasonable approximation. Note that the provided methodology assumes a one-to-one relationship between the performance metrics and the initial mass; in case there are multiple trajectory options for the same arc and the same initial mass (e.g., different coasting strategies), a sub-problem can be solved to identify one surrogate model to be used for the arc; or alternatively, multiple corresponding surrogate models can be considered as multi-graph arcs for the optimizer to make the choice.

\subsubsection{Mass Transformation}
The primary performance metric of any trajectory is measured in terms of deliverable mass fraction (i.e., final-to-initial mass ratio). Since this ratio itself may depend on the initial mass, a general relationship of the form $y^-_{ijve} = g(y^+_{ijve})$ is formed between the total final mass on an arc (i.e., inflow into the arc's destination node: $y^-_{ijve}$ ) and its total initial mass (i.e., outflow from the arc's origin node: $y^+_{ijve}$). Let us consider here how a linear mass transformation constraint (Equation~\eqref{eq:xformation1}) can be formulated for different forms of the function $g$.

\begin{enumerate}[label=(\alph*)]
\item Affine/linear function: First, consider the case where the relationship between the final and initial arc masses can be approximated as:
\begin{equation}
\label{eq:linear1}
y^-_{ijve} = g(y^+_{ijve}) = {}^1p\cdot y^+_{ijve} + {}^0p,
\end{equation}
where ${}^1p$ and ${}^0p$ are constants. Note that ${}^0p=0$ if the Tsiolkovsky rocket equation is satisfied. The above equation must be modified to reflect zero final mass when there is no initial mass so that it is consistent with physical reality. This is done by setting up the problem such that the component of the vector $\bm{x}^+_{ijve}$ that denotes the propulsive vehicle across the corresponding transportation arc is chosen to be binary-valued (or by introducing a new binary-valued variable that is dependent on it). Thus, if we denote this binary variable as $x^{\rm vehicle}_{ijve}$, the above relationship can be modified to:
\begin{equation}
\label{eq:linear2}
y^-_{ijve} = {}^1p\cdot y^+_{ijve} + {}^0p\cdot x^{\rm vehicle}_{ijve}.
\end{equation}
Consider an example with three commodities, where the vehicle is represented using integer variables in number of units, whereas the payload and propellant variables are represented using continuous variables measured in mass units. Assuming that the only commodity transformation that occurs on a transportation arc is the consumption of propellant, then the mass transformation constraint for the space logistics problem is implemented as:
\begin{equation} \footnotesize
\label{eq:xform_affine}
\left\lbrace\begin{bmatrix} 1 & 0 & 0 \\ 0 & 1 & 0 \\ {}^1p-1 & {}^1p-1 & {}^1p \end{bmatrix} \times \bm{M} + \begin{bmatrix} 0 & 0 & 0 \\ 0 & 0 & 0 \\ 0 & {}^0p & 0 \end{bmatrix} \right\rbrace\begin{bmatrix} \mathrm{payload} \\ \mathrm{vehicle} \\ \mathrm{propellant} \end{bmatrix}^+_{ijve}
 =
\bm{M} \times \begin{bmatrix} \mathrm{payload} \\ \mathrm{vehicle} \\ \mathrm{propellant} \end{bmatrix}^-_{ijve}. \normalsize
\end{equation}
where $\bm{M}$ is a diagonal matrix that convert the integer variables (e.g., vehicle) into continuous mass in kilograms. In this example, $\bm{M}$ is
\begin{equation} \footnotesize
\bm{M}=  \begin{bmatrix} 1 & 0 & 0 \\ 0 & m^{\rm vehicle} & 0 \\ 0 & 0 & 1 \end{bmatrix} 
\end{equation}
where $m^{\rm vehicle}$ is the dry mass of the vehicle (in mass units).

\item Nonlinear function: In a more general case where the relationship between the final mass and initial mass on an arc is written as a nonlinear function, the campaign logistics problem may be solved using mixed-integer nonlinear programming (MINLP) techniques; however such methods are usually computationally burdensome, not conducive to quick solutions, and often heuristics-based. Thus, to solve such problems in a computationally efficient manner within the event-driven network model, the MINLP structure is reduced to a MILP form using piecewise linear functions. 
\begin{figure}[htbp!]
\centering
\includegraphics[width=0.45\textwidth]{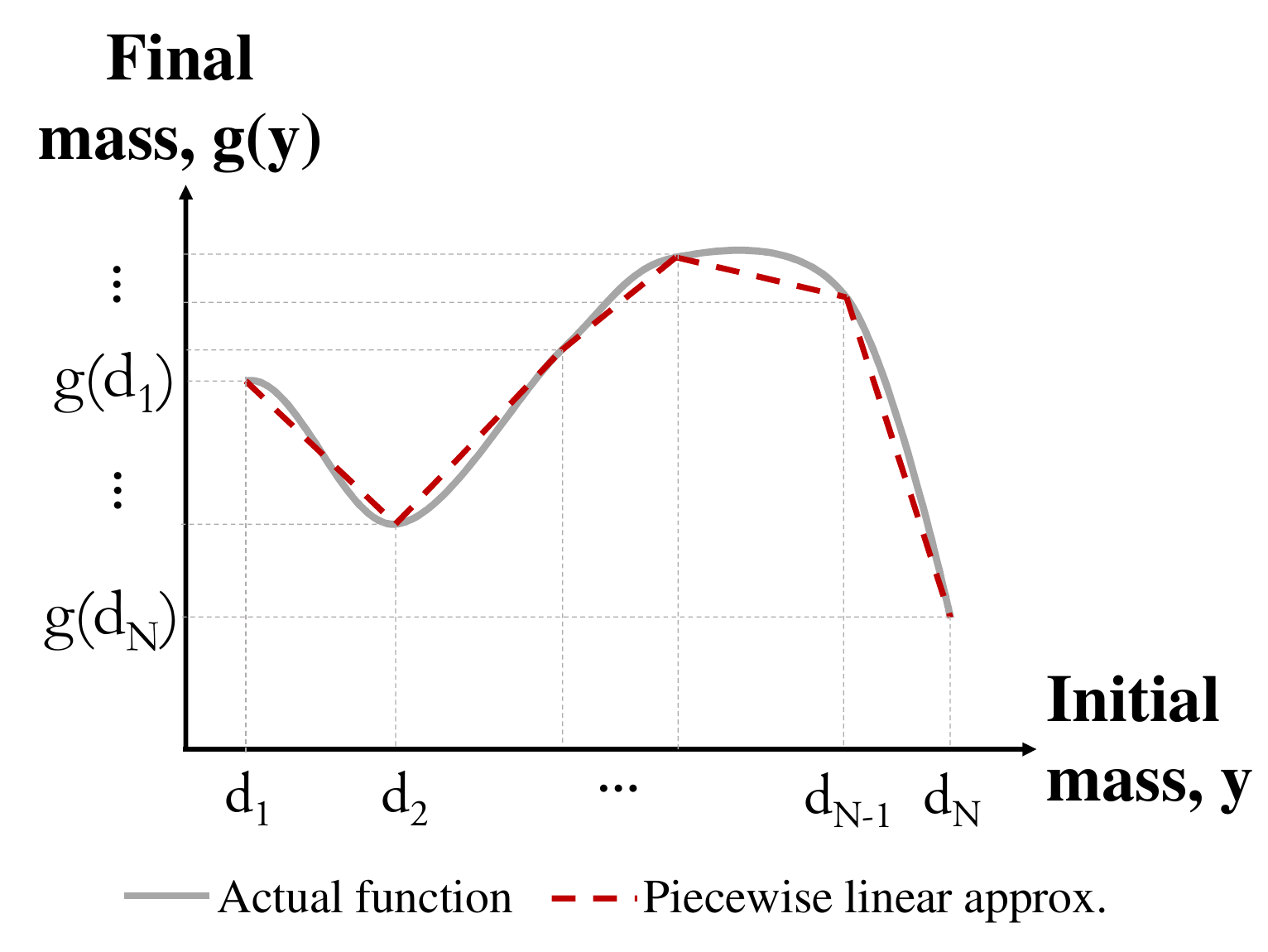}
\caption{Piecewise linear approximation of trajectory performance metrics vs. initial mass.\label{fig:pwl}}
\end{figure}

A general nonlinear mass transformation constraint $g(y)$ such as the one shown by the solid grey line in Figure~\ref{fig:pwl} can be approximated by piecewise linear (PWL) functions (red dashed lines) using the breakpoints \{$d_1, d_2$, ... $d_N$\}. If a set of additional continuous variables \{$\lambda_1, \lambda_2$, ... $\lambda_N$\} that belong to an SOS2 set\footnote{An SOS2 (special ordered sets of type 2) constraint enforces that at most two of the $\lambda$ can be nonzero and these two nonzero elements of the set must be consecutive.} are introduced, then the approximated value of the final mass $y^-_{ijve}$ at a desired initial mass of $y^+_{ijve}$ can be found as \cite{beale}:
\begin{subequations}
\label{eq:pwl}
\begin{align}
  \lambda_1 + ... + \lambda_N &= 1 \label{eq:pwl1}\hquad, \\
  \lambda_1d_1 + ... + \lambda_Nd_N &= y^+_{ijve} \label{eq:pwl2}\hquad, \\
  \lambda_1g(d_1) + ... + \lambda_Ng(d_N) &=  y^-_{ijve} \label{eq:pwl3}\hquad.
\end{align}
\end{subequations}
\noindent This set of equations can be used for each arc where a nonlinear performance metric model is encountered by adding the requisite number of $\lambda$ variables associated with every arc into the decision variables. The constraints in Equation~\eqref{eq:pwl} are also appended to the network formulation in Equations~\eqref{eq:cost}--\eqref{eq:tot_mass}. In addition, Equation~\eqref{eq:xformation1} is modified to conserve all commodities except the propellant:
\begin{equation}\footnotesize
\label{eq:xform_curve3}
\begin{bmatrix} 1 & 0 \\ 0 & 1  \end{bmatrix}_{ijv} \times \bm{M} \times \begin{bmatrix} \mathrm{payload} \\ \mathrm{vehicle} \end{bmatrix}^+_{ijve}
 =
\bm{M} \times \begin{bmatrix} \mathrm{payload} \\ \mathrm{vehicle} \end{bmatrix}^-_{ijve}. \normalsize
\end{equation}

Thus, by splitting the mass transformation constraint between Equation~\eqref{eq:pwl} and \eqref{eq:xform_curve3} in this case, the propellant consumption is implicitly considered even though the commodity denoting propellant is not explicitly shown.

\end{enumerate}

As an important note, the PWL approximation method described above in Equation~\eqref{eq:pwl} is one of the most basic models for linearizing MINLP problems. Other methods such as those reviewed in Reference~\cite{pwl} may also be employed as suited to the problem at hand. For example, PWL formulations developed by Vielma et al. \cite{vielma1,vielma2} are known to be more efficient in cases where a large number of breakpoints are considered. The accuracy of the PWL approximation depends on the number and location of the breakpoints. This breakpoint selection to minimize error is a separate optimization problem on its own and the interested reader is directed to studies such as \cite{pwl_pts1} and \cite{pwl_pts2} regarding this subject.

\subsubsection{Time of Flight }
The time required to transport commodities between two nodes is the secondary performance metric of a spacecraft trajectory. If the time of flight, or $\mathrm{TOF}$, across the trajectory that represents an arc, can be expressed as $\mathrm{TOF} = h(y^+_{ijve})$, then the arc length can be related to its commodity flow through this function $h$. Linear and nonlinear forms of $h$ can be converted to linear arc length relationships as discussed below. Once the arc lengths in the event-driven network are expressed as linear functions of the commodity flow, these relationships can be used to derive the campaign-specific time constraints that comprise Equation~\eqref{eq:tof_con}.

\begin{enumerate}[label=(\alph*)]
    \item Affine/linear function: First consider the case where the transportation arc can be traversed in a time span that varies linearly with the spacecraft's initial mass (i.e. total outflow mass from the arc's origin node) as
    \begin{equation}
        \mathrm{TOF} = h(y^+_{ijve}) = {}^1q\cdot y^+_{ijve} + {}^0q,
    \end{equation}
   \noindent where ${}^1q$ and ${}^0q$ are both constants derived from trajectory analysis. As with the mass transformation constraint, this is modified to reflect zero arc length for no initial mass:
    \begin{equation}
        \Delta t_{ijve} = {}^1q\cdot y^+_{ijve} + {}^0q\cdot x^{\rm vehicle}_{ijve},
        \label{eq:arcl}
    \end{equation}
    \noindent where $\Delta t_{ijve}$ is the arc length and $x^{\rm vehicle}_{ijve}$ is the binary-valued component of the commodities vector $\bm{x}^+_{ijve}$ specifying the flow of vehicle $v$ over the arc ($i,j,v,e$).
    
    \item Nonlinear function: Conversely, if $\mathrm{TOF}$ is some general nonlinear function of the initial spacecraft mass, then PWL functions can be utilized again to derive the arc length in a manner similar to that shown earlier in Equation~\eqref{eq:pwl}. If the nonlinear curve $h$ is approximated by piecewise linear segments using breakpoints \{$d_1, d_2, ... d_N$\} and additional continuous SOS2 variables \{$\lambda'_1, \lambda'_2, ... \lambda'_N$\} are introduced, then the arc length $\Delta t_{ijve}$ at a desired initial mass of $y^+_{ijve}$ can be approximated as:
    \begin{subequations}
    \label{eq:pwl_t}
    \begin{align}
      \lambda'_1 + ... + \lambda'_N &= 1 \label{eq:pwl1_t}\hquad, \\
      \lambda'_1d_1 + ... + \lambda'_Nd_N &= y^+_{ijve} \label{eq:pwl2_t}\hquad, \\
      \lambda'_1h(d_1) + ... + \lambda'_Nh(d_N) &= \Delta t_{ijve} \label{eq:pwl3_t}\hquad.
    \end{align}
    \end{subequations}
    Equation~\eqref{eq:pwl_t} should then be appended to the original event-driven space logistics formulation. The set of $\lambda'$ variables should also also be added to the decision variables.
\end{enumerate}

This completes the derivation of the MILP formulation for the event-driven space logistics optimization framework, where variable flow-dependent arc lengths are incorporated by using event-based time steps. This model is derived based on the assumption that mass transformations occurring over the holdover arcs do not depend on time. Modification to the proposed method is needed if time-dependent transformations over holdover arcs, such as ISRU, are to be considered. The rest of this paper presents and solves a case study using the newly developed event-driven network formulation, thereby demonstrating the use of the campaign-specific time-related constraints.

\section{Case Study Overview \label{sec:case}}
NASA's journey to Mars is expected to go through near-term missions to the cislunar space, which will allow human space programs to progressively fill knowledge gaps and prove newly developed technologies. Within this context, propellant resupply depots have been considered a popular logistics strategy for enabling ambitious manned missions beyond Low Earth Orbit \cite{depot1,depot2,depot3,depot4,depot_ho}. Furthermore, previous architectural concepts have proposed using the collinear Earth–moon Lagrange points to stage mission infrastructure elements \cite{depot_ho,basecamp}. Many researchers have also concentrated on designing optimal high-thrust \cite{howell,mingtao,shah} and low-thrust trajectories within the cislunar system \cite{mingotti,martin}. Some studies have also attempted to combine low-thrust propulsion with high-thrust chemical propulsion, albeit for a single leg of the mission, and not for a campaign \cite{kluever1, kluever2, kluever3}. Some of these referenced studies suggest ideal locations for refueling the crew vehicle based on intuition or manual analysis, while other analyses optimize mission costs and design trajectories with respect to an assumed architecture. Thus, existing frameworks can lead to suboptimal design solutions without exploring the combined tradespace of the campaign architecture and trajectory/propulsion options. The choice of the case study campaign to demonstrate the newly developed space logistics framework developed in this paper is driven by this background. The campaign considered in this paper is the support of crewed Apollo-like missions to the lunar surface by using in-space refueling tanks pre-deployed by a fleet of high-thrust chemical propulsion (CP) and low-thrust solar electric propulsion (SEP) tugs. This section describes in detail the problem setup for this case study campaign. 

\subsection{Problem Description}
First, consider a series of repeating Apollo-style crewed missions for lunar surface exploration utilizing only the carry-along logistics strategy. For each mission, the crew stack (comprised of the Crew \& Service Module, or CSM, and the Lunar Module, or LM) departs low Earth orbit via a translunar injection maneuver provided by the Saturn V's upper stage. The CSM inserts the entire stack into lunar orbit once in the vicinity of the moon. Lunar ascent and descent operations are then conducted by the LM, which is discarded once the crew is returned to the CSM idling in low lunar orbit. Finally, the crew is returned to Earth by the trans-Earth injection maneuver provided by the CSM. In this scenario, the crew stack carries its entire propellant supply. 

Instead, consider a campaign composed of two stages -- the cargo delivery phase, which then facilitates the crew mission phase. In the cargo delivery phase, a fleet of cargo tugs delivers the crew resupply propellant to way stations in the form of disposable ``droptanks''. These refueling stations are delivered using propulsion or trajectory options that are not normally considered for crewed spacecraft, e.g., using cheaper transfers that may be long duration or using SEP tugs. 
Once the first phase is complete, the crew mission/s can utilize these in-space propellant stockpiles in the second phase. Using this distributed-launch architecture lets mission designers circumvent the necessity of using super heavy-lift launch vehicles for crew launch. This logistics problem can be optimized along multiple dimensions, three of which are 1) locations of the way stations where the crewed spacecraft can be refueled, 2) amount of propellant resupplied to the crew at each possible way station, and 3) the tugs and their propulsion system (i.e. high-thrust or low-thrust) selected for delivery of these propellant resupply tanks. The event-driven space logistics model developed in this paper can solve this problem along all three of these dimensions.

\begin{figure}[htb!]
\centering
\includegraphics[width=0.9\textwidth]{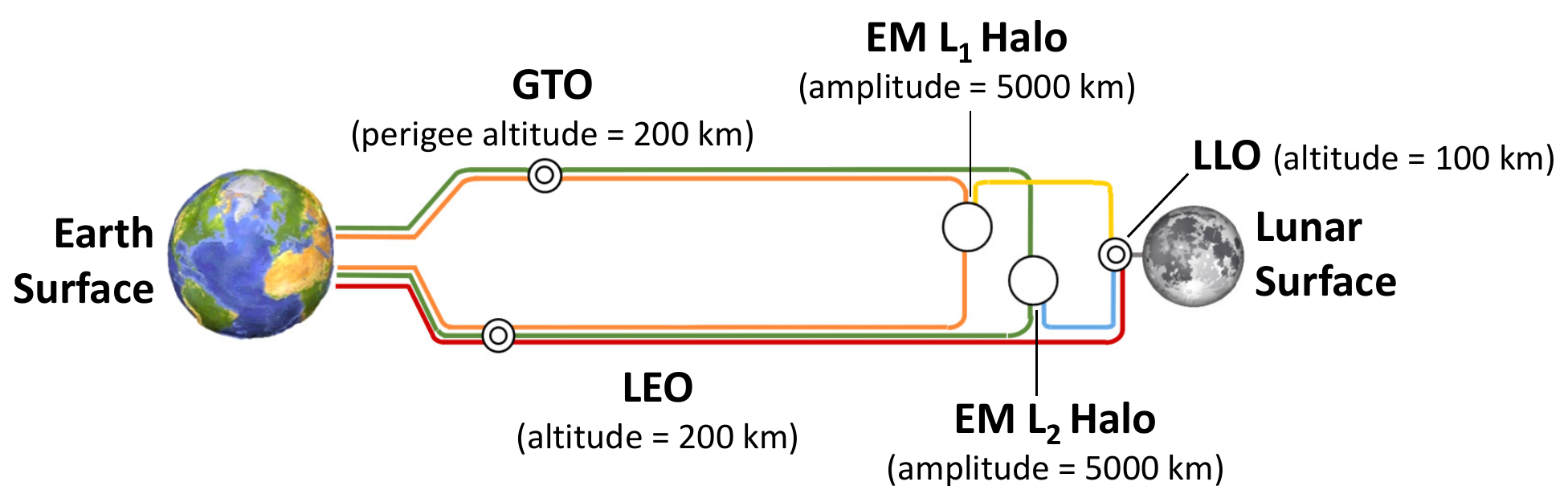}
\caption{Network considered for case study campaign.\label{fig:static}}
\end{figure}

Figure~\ref{fig:static} summarizes the static network considered for this cislunar case study problem. In the current paper, halo orbits at the Earth-moon Lagrange points EML$_1$ and EML$_2$ are considered as possible way stations for refueling the crew stack as they have associated invariant manifolds that can be used to provide cheaper low-energy pathways (i.e. lower $\Delta$v, but longer time of flight) for the tug journeys. The cargo delivery fleet is considered to only consist of discrete vehicle choices as specified in Table~\ref{tab:tugs}, with propulsion systems of high-thrust (chemical) or low-thrust (solar electric) type. This means that the dry mass of each tug and its fuel capacity is fixed and the payload mass carried by a tug across the network is only limited by its fuel capacity. Additionally, multiple tugs of the same kind are allowed in the analysis. Each tug is also permitted to be reused twice (i.e., a total of three uses per tug), so that the number of launches required can be minimal. In order to be reused, a tug must return to its corresponding Earth parking orbit to be refueled. In such cases, the tug would be required to carry additional fuel for facilitating its return, in addition to the cargo it needs to deliver. 

\begin{table}[h]
\centering
\caption{Specifications of tugs used in case study campaign. \label{tab:tugs}}
\begin{threeparttable}
\begin{tabular}{l>{\centering\arraybackslash}p{16mm}>{\centering\arraybackslash}p{16mm}>{\centering\arraybackslash}p{16mm}c>{\centering\arraybackslash}p{16mm}>{\centering\arraybackslash}p{16mm}>{\centering\arraybackslash}p{16mm}}
\hline \hline
 & \multicolumn{3}{c}{Chemical Propulsion (CP) tugs} & & \multicolumn{3}{c}{Solar Electric Propulsion (SEP) tugs} \\
 \cline{2-4} \cline{6-8}
 & Type 1 \cite{tug1}& Type 2 \cite{tug2}& Type 3 \cite{tug3}& & Type 4 \cite{tug4}& Type 5 \cite{tug5}& Type 6 \cite{tug6} \\
 \hline
Dry mass, t & 2.3 & 5.5 & 6.5 && 3.5 & 7.68 & 10.7 \\
Propellant capacity, t & 11.5 & 41 & 68 && 11 & 16 & 39 \\
Power, kW\tnote{+} & - & - & - && 40 & 150 & 300 \\
$I_\text{ sp}$, s & 450 & 450 & 450 && 3000 & 3000 & 2000 \\
Number of units & 2 & 2 & 3 && 2 & 2 & 1 \\
\hline \hline
\end{tabular}  
  \begin{tablenotes}
    \item[+] \small Efficiency in power transfer to propulsion system assumed to be 60\% to calculate thrust provided by SEP tugs. 
  \end{tablenotes}
\end{threeparttable}
\end{table}

The in-space propellant resupply chain for the crew spacecraft is modeled in the form of disposable droptanks that are sized by their fuel storage requirements through the use of the fuel's structural coefficient ($\varepsilon$). In the current analysis, the cost of rendezvous between the crew spacecraft and the droptanks at way stations is ignored. Stationkeeping costs for maintaining the droptanks at their operational locations (i.e. at the halo orbits or at the low-lunar orbit (LLO)) are also disregarded. Though not modeled here, the crew can be resupplied with other necessary commodities (such as life-support supplies) at these waystations; thus, the extended crew flight times due to the need to transit through way stations for the purposes of refueling is assumed to not have any appreciable effect on the system-level design of the crew spacecraft. The crew vehicles used in this analysis -- the CSM, LM and the launch vehicle upper stage -- are derived from Apollo missions' data. Note that the crew vehicles use high-density storable propellants and hence the in-space storage of refuel propellant is modeled using simple droptanks. Furthermore, the upper stage of the crew launch vehicle is assumed to be sized according to the size of the crew stack. 
Given these modeling parameters, the list of commodities required to simulate this case study campaign is provided in Table~\ref{tab:commodities}. 

\begin{table}[htb!]
\centering
\caption{List of commodities. \label{tab:commodities}}
\begin{tabular}{rlll}
\hline \hline
Commodity name & Abbreviation & Variable type & Notes\\
\hline
Upper stage structure & \code{strUS} & Continuous & - \\
Upper stage fuel & \code{fUS} & Continuous & - \\
CSM & \code{CSM}  & Binary & - \\
CSM fuel & \code{fCSM} & Continuous & - \\
LM & \code{LM} & Binary & - \\
LM fuel &\code{fLM} & Continuous & - \\
Droptank structure & \code{strDtank} & Continuous & - \\
CP tug fuel & \code{fHIGH} & Continuous & -\\
SEP tug fuel & \code{fLOW} & Continuous & - \\
Tug \#1 & \code{tug1} & Binary & Chemical tug type 1 \\
Tug \#2 & \code{tug2} & Binary & Chemical tug type 1 \\
Tug \#3 & \code{tug3} & Binary & Chemical tug type 2 \\
Tug \#4 & \code{tug4} & Binary & Chemical tug type 2 \\
Tug \#5 & \code{tug5} & Binary & Chemical tug type 3 \\
Tug \#6 & \code{tug6} & Binary & Chemical tug type 3 \\
Tug \#7 & \code{tug7} & Binary & Chemical tug type 3 \\
Tug \#8 & \code{tug8} & Binary & SEP tug type 1 \\
Tug \#9 & \code{tug9} & Binary & SEP tug type 1 \\
Tug \#10 & \code{tug10} & Binary & SEP tug type 2 \\
Tug \#11 & \code{tug11} & Binary & SEP tug type 2 \\
Tug \#12 & \code{tug12} & Binary & SEP tug type 3 \\
\hline \hline
\end{tabular}
\end{table}

The total cost of the campaign is measured in terms of the initial mass to low-Earth orbit (IMLEO). In the current problem context, an additional cost arises in the form of each mission's operational time. This indirect cost regarding the campaign duration has two components -- the total time required to set up the in-space refuel supply chain using the cargo delivery tugs and the total crew flight time across all considered human missions. These two time-related pseudo-costs are included in the event-driven network optimization framework as constraints (instead of as direct costs), so that the tradeoff between IMLEO and campaign duration can be explored. Deriving and optimizing the mathematical model of this problem can help mission designers explore combinations of logistics strategies of carry-along and refuel, while optimally combining propulsion technologies within these strategies. This mathematical model needs to be supplied with the trajectory performance surrogate models for the campaign, which are described next.

\subsection{Trajectory Performance Surrogate Models}
The performance of transporting commodities over an arc depends on the functions that relate the initial mass to the final mass and time of flight, and these functions depend on the propulsive element used on the arc and the trajectory followed. All trajectories in this campaign are considered within the circular restricted three-body model and some of them make use of the dynamical structures available in the Earth-moon system. The interaction of the gravitational pulls of the Earth and the moon creates families of dynamically unstable periodic orbits around the system's Lagrange points, giving rise to the invariant manifolds that act as low-energy pathways. These cheap transport options typically take longer than direct trajectories between the Earth and the moon, but can be especially useful for delivering cargo. Thus, all cargo deliveries use the low-energy pathways associated with the chosen halo orbits. As a consequence, the tugs cannot access LLO directly from their Earth parking orbits and must instead transfer through either the EML$_1$ or EML$_2$ halo orbit. On the other hand, crew transport is conducted using only high-thrust chemical propulsion on direct pathways that do not use the low-energy manifolds. 

In the current study, the performance of high-thrust transfers is taken directly from available literature. However, the low-thrust propulsion performance surrogate model is derived for this case study and supplied to the logistics framework. The contributions to the total cost of the campaign, divided up between the crew transport and the cargo transport, are discussed next.

\subsubsection{Crew Transfer} 
Past studies have designed minimum $\Delta$v trajectories suitable for human missions from low-Earth orbit (LEO) to LLO, as well as to periodic orbits around the EML$_1$ and EML$_2$ \cite{man2eml1,man2eml2,man2eml3,man2eml4}. In order to restrict the crew flight time to under safe limits, these transfers do not exploit the low-energy pathways, instead choosing to directly transfer onto the target orbits or to use lunar flybys. These crew-specific spaceflight routes also generally leverage Apollo-like free-return trajectories in order to enable safe aborts. The $\Delta$v values for the different arcs that can be traversed by the crew vehicles are gathered from available literature and summarized in Table~\ref{tab:crewcosts} of Appendix A. 

While the crew vehicles that provide transport in space have fixed structural masses (in terms of dry mass and fuel capacity), the upper stage of the crew launch vehicle is sized according to the fuel needed to provide the impulse demanded by the crew stack, i.e., commensurate to the magnitude of the translunar injection maneuver (TLI). This TLI varies depending on whether the campaign-level architectural solution decided by the logistics framework dictates the first stop to be at LLO, at the EML$_1$ halo orbit, or at the EML$_2$ halo orbit. Once this first impulse is completed by the upper stage, all succeeding maneuvers are conducted by the CSM or the LM (for lunar landing/ascent operations only). Hence, the crew only needs to be resupplied either with CSM fuel or LM fuel at any in-space way station. The details regarding these crew vehicles are provided in Table~\ref{tab:crewveh} of Appendix B. 

\subsubsection{Cargo Delivery for High-Thrust CP Tugs}
Apart from the crew transportation trajectory, the other contribution to the total IMLEO cost arises from the launch of tugs and their corresponding cargo. Any resupply propellant (i.e., CSM or LM fuel) stored in space is assumed to be held in disposable droptanks, which are sized linearly with the amount of propellant they are required to store. The structural mass of these in-space droptanks is calculated using the respective fuel's structural coefficient. When delivering these droptanks to their in-space storage locations, the high-thrust chemical propulsion tugs use the cheaper (but slower) routes through the Earth-moon system. The $\Delta$v data for cargo deliveries conducted using CP tugs is provided in Table \ref{tab:cptug} (Appendix A). These high-thrust tugs are assumed to begin their journey in LEO. All tugs can be reused during the campaign if their own propellant tanks are refilled; however, they must return to their Earth parking orbit (which is LEO for the case of CP tugs) to do so. 

\subsubsection{Cargo Delivery for Low-Thrust SEP Tugs}
The propellant mass consumed by an SEP spacecraft and its time of flight are heavily dependent on the ratio of the thrust it can provide to its total initial mass. This is because low-thrust propulsion consists of continuous burn maneuvers (as opposed to impulsive burn approximation used for chemical rocket maneuvers). Thus, the performance of low-thrust transfers cannot be simply be collated from existing literature. Obtaining these low-thrust trajectories is a computationally intensive process involving design of the many-revolution thruster-on trajectory segment required to depart/arrive at the Earth/moon \cite{martin,ozimek2010,mingotti2007,mingotti2009,mingotti2011,mingotti2012}. 
Since most existing low-thrust trajectory design methods did not lend themselves well to quick parametric trade studies, our previous studies developed a computationally efficient and reasonably accurate method for calculating low-thrust transportation performance metrics over the considered arcs in the cislunar region \cite{jag_c2,jag_j2}. The interested reader is directed to these previous works for a complete description, which is summarized here for clarity and continuity. The trajectories used by the low-thrust cargo arcs are composed of two connected segments -- the continuous constant-thrust spiral segment, and the coast segment on the invariant manifold. The multi-revolution (spiral) portion of the trajectory assumes continuous constant thrust and constant $I_\text{sp}$ as is common in the previously cited literature that addresses the same trajectory design problem. This trajectory segment is treated as an orbital transfer and thus designed using a feedback control law called Q-law \cite{petro2004,petro2005}. Feedback control laws are known to be suboptimal, but they simplify the solution process and provide good estimates for low-thrust transfer performance metrics, thus making the method favorable to large-scale trade studies or for quick approximations such as the application here. In order to use Q-law with dynamical systems theory to obtain low-thrust arc performance metrics, favorable points on the manifold must be first chosen to connect it with the spiral trajectory segment. These manifold patch points, which separate the continuous-thrust spiral arc from the coast arc, are chosen via a grid search across three different values of thrust accelerations \cite{jag_c2}. The Geosynchronous Transfer Orbit (GTO) is chosen as the initial Earth parking orbit for the low-thrust tugs, in compliance with the selection made by other previous studies that studied similar trajectories \cite{mingotti2007,martin}. 
Once the manifold patch points are chosen, the trajectories for various combinations of low-thrust tugs from Table~\ref{tab:tugs} and their payloads are designed for all permissible arcs. The final mass and time of flight are calculated for each SEP tug, by starting at the tug's dry mass as the initial mass along the arc and increasing its payload mass in discrete steps up to the maximum value permissible by the tug's fuel capacity. The trajectory performance models for each low-thrust arc in the network are shown in Figure~\ref{fig:fit} \footnote{See Appendix C for discussion}. 

\begin{figure}[ht]
\centering
\includegraphics[width=0.55\textwidth]{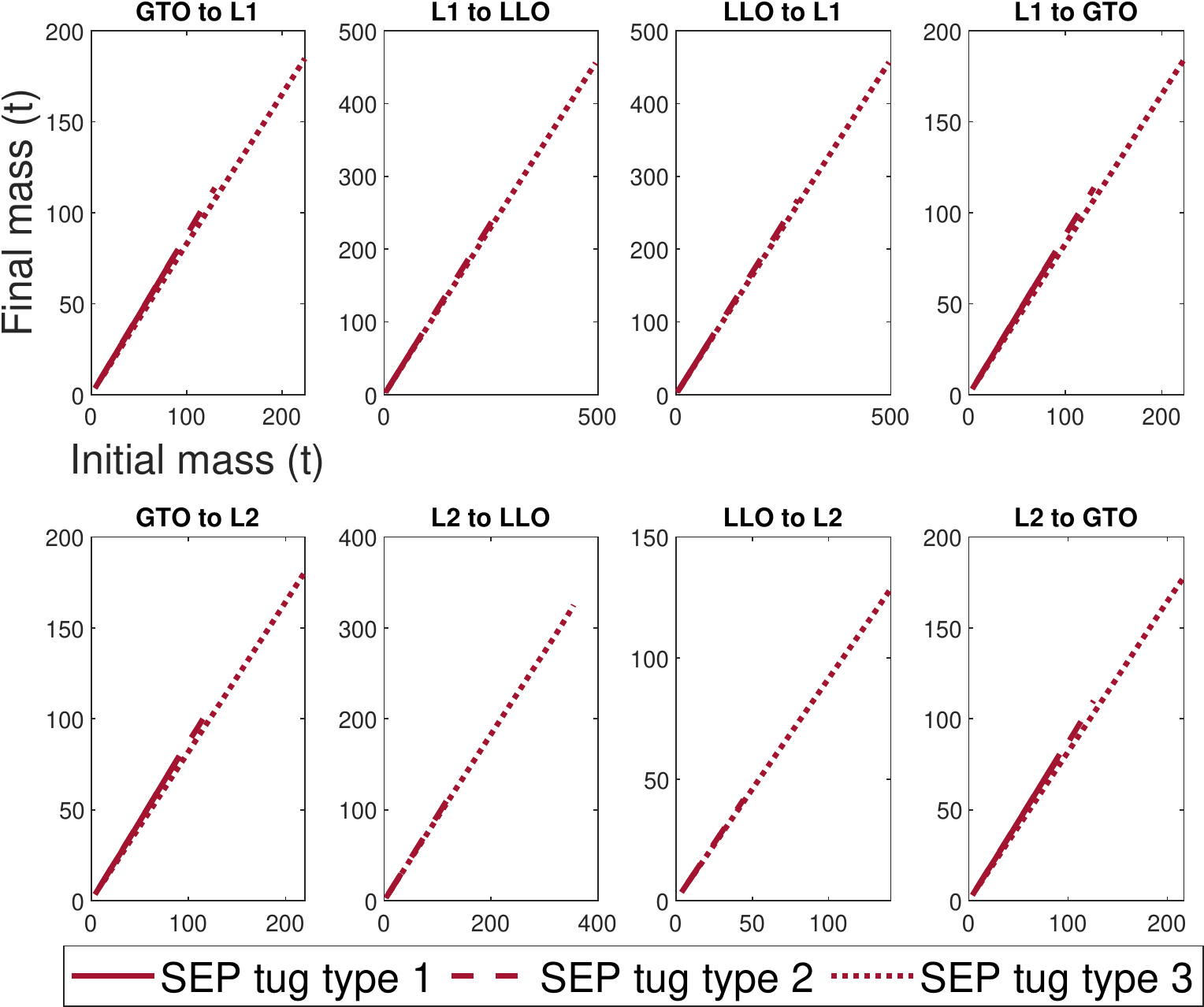}
\includegraphics[width=0.55\textwidth]{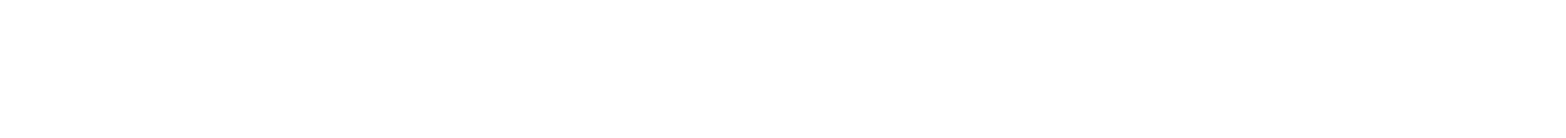}
\includegraphics[width=0.55\textwidth]{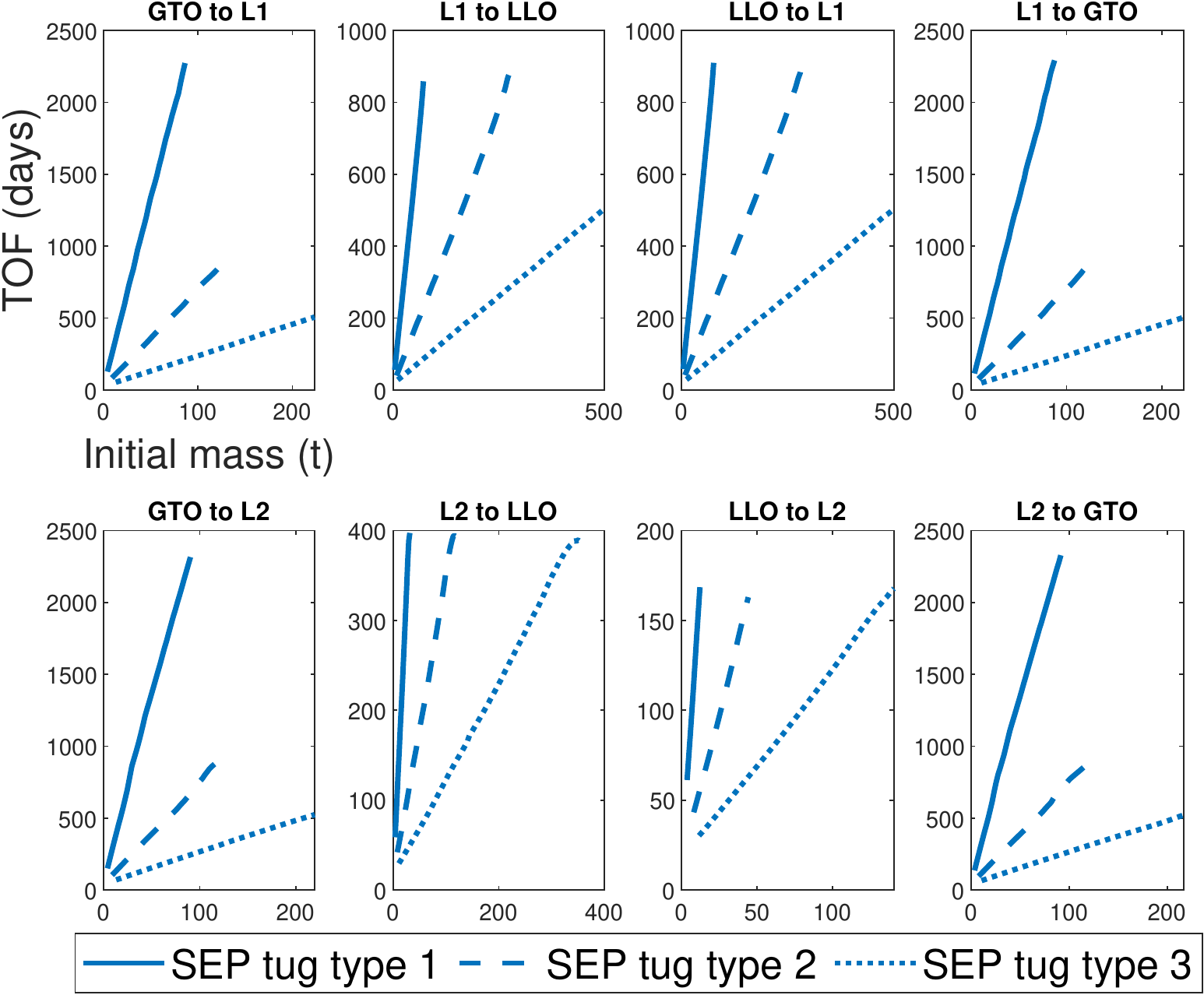}
\caption{Final mass (in t) and time of flight (TOF; in days) vs. initial mass (in t) for each SEP tug in Table \ref{tab:tugs} \label{fig:fit}}
\end{figure}

Note that although the above models are used as an example case study in this paper, other models obtained through a more-detailed low-thrust trajectory evaluation can also be incorporated in the proposed event-driven space logistics optimization framework. 
For example, higher mass savings can be obtained if coast arcs are allowed in the spiral portion of the low-thrust cislunar trajectories. Additionally, alternative initial orbits besides GTO can also be explored. 

\subsubsection{Linear Fit for the Trajectory Performance Models}
In the case study, high-thrust transfers use constant $\Delta$v (which indicates constant final-to-initial mass ratio for each vehicle) and constant time of flight for each arc, whereas low-thrust transfers use a linear approximation model to relate the initial mass to the final mass and time of flight (refer to Figure~\ref{fig:fit}). Thus, a general linear relationship can be considered:
\vspace{-5mm}
\begin{align}
\label{eq:massfit}
\newsubeqblock
\mysubeq y^-_{ijve}= {}^1p\cdot y^+_{ijve} + {}^0p \\
\mysubeq \mathrm{TOF}= {}^1q\cdot y^+_{ijve} + {}^0q 
\label{eq:timefit}
\end{align}

\noindent where $p$ and $q$ are constant coefficients which are derived from trajectory analysis and used later to incorporate these performance models into the event-driven space logistics model. For high-thrust transfers, ${}^0p={}^1q=0$ and (${}^1p,{}^0q$) take the corresponding values from Table~\ref{tab:cptug}. For low-thrust transfers permitted in the case study network, the coefficients (${}^1p,{}^0p,{}^1q,{}^0q$) are derived from fitting a linear curve to the performance metrics derived in Figure~\ref{fig:fit} and are listed in Tables~\ref{tab:fit_mass}--\ref{tab:fit_tof}. These relationships are modified into the form of Equations~\eqref{eq:linear2} and \eqref{eq:arcl} when implemented as constraints in the MILP formulation, as explained later in this section. As before, the initial spacecraft mass on the arc (which is equivalent to the total outflow from the arc's origin node) is denoted by $y^+_{ijve}$ and the final spacecraft mass (which is equivalent to the total inflow into the arc's destination node) is represented by $y^-_{ijve}$.

\begin{table}[htb!]
\centering
\caption{Low-thrust transfer performance metrics: Linear fit coefficients for final mass (in t) vs. initial mass (in t) \label{tab:fit_mass}}
\begin{tabular}{rccccccccc}
\hline \hline
 Arc & & \multicolumn{2}{c}{SEP tug type 1} & & \multicolumn{2}{c}{SEP tug type 2} & & \multicolumn{2}{c}{SEP tug type 3} \\
 & & ${}^1p$ & ${}^0p$ & & ${}^1p$ & ${}^0p$ & & ${}^1p$ & ${}^0p$ \\
\cline{1-1} \cline{3-4} \cline{6-7} \cline{9-10}
\\[-8pt]
\multicolumn{10}{c}{\textit{To/\ from EML$_1$ halo orbit}} \\
GTO to L1 & & 0.8757 & -0.0038 & & 0.8772 & 0.1202 & & 0.8251 & 0.2401 \\
L1 to LLO & & 0.9446 & 0.0446 & & 0.9447 & 0.1586 & & 0.9214 & -0.0391 \\
LLO to L1 & & 0.9446 & 0.0429 & & 0.9449 & 0.1398 & & 0.9214 & -0.0407 \\
L1 to GTO & & 0.8758 & -0.0082 & & 0.8764 & 0.1453 & & 0.8253 & 0.1808 \\[3pt]
\multicolumn{10}{c}{\textit{To/\ from EML$_2$ halo orbit}} \\
GTO to L2 & & 0.8801 & -0.2208 & & 0.8709 & 0.1763 & & 0.8219 & -0.2458 \\
L2 to LLO & & 0.9397 & 0.0338 & & 0.9396 & 0.1337 & & 0.9112 & 0.5247 \\
LLO to L2 & & 0.9394 & 0.0276 & & 0.9402 & 0.0841 & & 0.9128 & 0.2931 \\
L2 to GTO & & 0.8818 & -0.2571 & & 0.8718 & 0.1172 & & 0.8232 & -0.3195 \\
\hline \hline
\end{tabular}
\end{table}

\begin{table}[htb!]
\centering
\caption{Low-thrust transfer performance metrics: Linear fit coefficients for time of flight (in days) vs. initial mass (in t) \label{tab:fit_tof}}
\begin{tabular}{rccccccccc}
\hline \hline
 Arc & & \multicolumn{2}{c}{SEP tug type 1} & & \multicolumn{2}{c}{SEP tug type 2} & & \multicolumn{2}{c}{SEP tug type 3} \\
 & & ${}^1q\ (\times 10$) & ${}^0q$ & & ${}^1q\ $ & ${}^0q$ & & ${}^1q\ $ & ${}^0q$ \\
\cline{1-1} \cline{3-4} \cline{6-7} \cline{9-10}
\\[-8pt]
\multicolumn{10}{c}{\textit{To/\ from EML$_1$ halo orbit}} \\
GTO to L1 & & 2.598 & 26.631 & & 6.832 & 19.146 & & 2.164 & 22.865 \\
L1 to LLO & & 1.156 & 8.567 & & 3.074 & 8.403 & & 0.973 & 17.465 \\
LLO to L1 & & 1.157 & 8.015 & & 3.062 & 9.199 & & 0.973 & 17.482 \\
L1 to GTO & & 2.595 & 27.177 & & 6.872 & 18.159 & & 2.162 & 23.674 \\[3pt]
\multicolumn{10}{c}{\textit{To/\ from EML$_2$ halo orbit}} \\
GTO to L2 & & 2.505 & 88.126 & & 7.180 & 32.179 & & 2.204 & 45.036 \\
L2 to LLO & & 1.257 & 13.027 & & 3.359 & 12.346 & & 1.098 & 13.092 \\
LLO to L2 & & 1.266 & 13.364 & & 3.329 & 14.458 & & 1.079 & 15.511 \\
L2 to GTO & & 2.476 & 91.578 & & 7.128 & 35.899 & & 2.188 & 45.804 \\
\hline \hline
\end{tabular}
\end{table}
\FloatBarrier

As a final point of note, one can observe that the initial Earth parking orbit for the low-thrust SEP tugs is GTO, while the campaign cost is measured in terms of IMLEO. To account for this discrepancy, any mass launched to GTO is penalized to be 1.74 times the initial mass in LEO; this penalty is derived from the ratio of initial mass to final mass if a rocket with $I_\text{sp}=450$ s is used to place a spacecraft in GTO when starting at LEO, which requires $\Delta\mathrm{v}\sim2.4554$ km/s.

\subsection{Event-Based Time Steps for Dynamic Network Modeling}
\begin{figure}[htb!]
\centering
\includegraphics[width=0.8\textwidth]{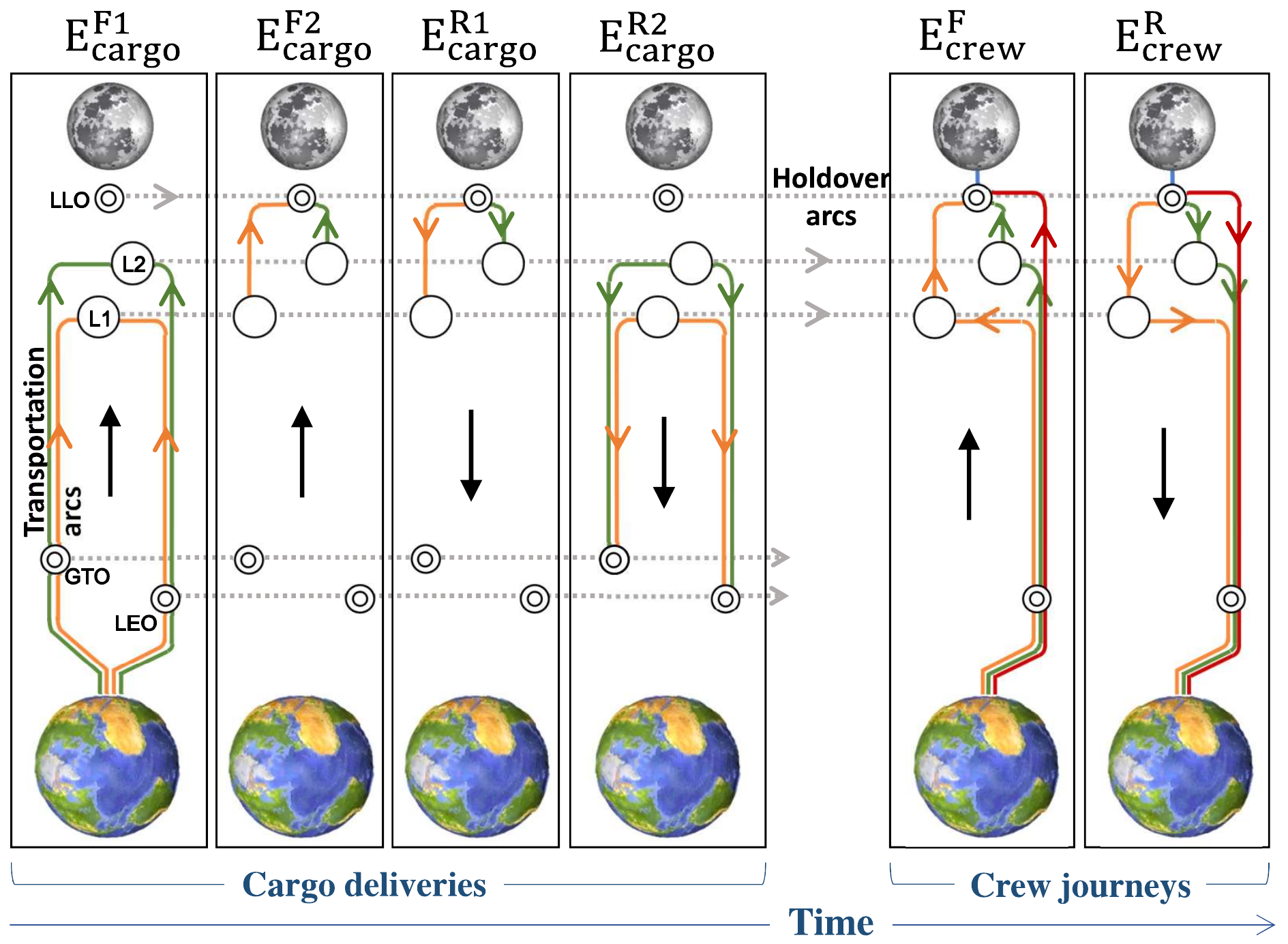}
\caption{Event-driven network for the case study campaign.\label{fig:edgmcnf}}
\end{figure}

As explained earlier, the static space logistics design model lacks any mechanism for dealing with temporal network behavior and can result in erroneous flow generation loops (among other problems). As opposed to the existing dynamic variant (i.e., the time-expanded network) where discrete time steps are used to expand the static network, event-steps are used in the current work to allow easy integration of low-thrust transfer times (which are not known \textit{a priori} and thus preclude the use of discrete time step sizes). Figure~\ref{fig:edgmcnf} shows the breakdown of events for the current case study, where four event-steps are used for each use of the cargo tugs and two additional event-steps for every crew mission. Each of these layers represents a copy of the complete static network -- they span a variable time period that is calculated and adjusted internally to meet the overall time constraints on the campaign duration. 

The first four static layers or event steps ($\rm{E^{F1}_{cargo}}$, $\rm{E^{F2}_{cargo}}$, $\rm{E^{R1}_{cargo}}$, and $\rm{E^{R2}_{cargo}}$) permit only cargo deliveries by tugs chosen from the available fleet. Layers $\rm{E^{F1-F2}_{cargo}}$ allow only those arcs that move toward the moon (``forward''), while $\rm{E^{R1-R2}_{cargo}}$ only allow returning arcs. Restricting transport directions in this way avoids flow generation loops. Parallel arcs (e.g., LEO $\rightarrow$ EML$_1$ halo and LEO $\rightarrow$ EML$_2$ halo) that do not interact with each other are allowed within the same layer. By allowing a single layer to encompass all flow until the end of an event (e.g., forward motion towards the halo orbits, in layer $\rm{E^{F1}_{cargo}}$), commodity interaction at nodes can be considered between the events. For example, if two tugs finish at the same node at the end of an event layer, then they could exchange cargo, and more importantly, they could even exchange propellant if both tugs employ the same type of propulsion systems (high-thrust or low-thrust). Thus, a demarcation between $\rm{E^{F1}_{cargo}}$ and $\rm{E^{F2}_{cargo}}$ layers is made so that multiple vehicles can interact with each other at an intermediate node (i.e., the EML$_1$ halo or the EML$_2$ halo). 

Similarly, the crew journeys are also split into layers $\mathrm{E^F_{\rm crew}}$ and $\mathrm{E^R_{\rm crew}}$. The four layers required for the tug flight are collapsed to two layers for crew travel because, as shown later, only a single crew vehicle path is considered, and thus there is no commodity exchange among the vehicles within these event layers.

The layers $\mathrm{E^{F1-F2}_{\rm cargo}}$ and $\mathrm{E^{R1-R2}_{\rm cargo}}$ are repeated three times in the current problem to correspond to the three uses allowed for each tug. Likewise, the layers $\mathrm{E^F_{\rm crew}}$ and $\mathrm{E^R_{\rm crew}}$ are repeated three times to signify the three crew missions that the current campaign is set up for. If $\mathbb{E}_{\rm cargo}$ consolidates all the cargo-related layers such that $\mathbb{E}_{\rm cargo} = \mathrm{E^{F1}_{\rm cargo}} \cup \mathrm{E^{F2}_{\rm cargo}} \cup \mathrm{E^{R1}_{\rm cargo}} \cup \mathrm{E^{R2}_{\rm cargo}}$, then the only propulsive elements allowed in the event set $\mathbb{E}_{\rm cargo}$ are the cargo tugs. Along the same lines, only crew vehicles are permitted to provide propulsion along the arcs in the event set $\mathbb{E}_{\rm crew} = \mathrm{E^F_{\rm crew}}\cup\mathrm{E^R_{\rm crew}}$. The cargo (i.e. refuel droptanks) delivery-related events are assumed to occur before the crew-related events and these two sets of event layers are connected via the holdover arcs. Together these two sets ($\mathbb{E}_{\rm cargo}$ and $\mathbb{E}_{\rm crew}$) form the set of all events used to discretize the time dimension, i.e., $\mathbb{E}= \mathbb{E_{\rm cargo}} \bigcup \mathbb{E_{\rm crew}}$.

Finally, the vehicles are also divided into two sets -- the set of cargo tugs $\mathbb{V}_\mathrm{cargo}$ and the set of crew vehicles $\mathbb{V}_\mathrm{crew}$ together form the set representing the entire vehicle fleet $\mathbb{V} = \mathbb{V}_\mathrm{cargo} \bigcup \mathbb{V}_\mathrm{crew}$. The flow of cargo tugs across the network is regulated through the use of multi-graphs (refer to Figure~\ref{fig:arc2}). Multiple transportation arcs connect a single pair of nodes within any event layer in the set $\mathbb{E}_{\rm cargo}$ to form the multi-graph between those nodes. Each arc in the multi-graph set is then assigned to a single tug unit, with the implication that the propulsion required to traverse a single arc among the multi-graph arcs between the same pair of nodes can only be provided by a particular tug unit. This is an important modeling feature that aids in the correct consolidation of flight times of different tugs across event layers and is further explained in subsequent sections.

\subsection{Demand and Supply}
The Earth's surface supplies a practically infinite amount of all other commodities that are continuous variables (such as the droptank structural mass, upper stage structural mass, and all propellants). Cost coefficients are thus only applied to those arcs that launch commodities from the Earth surface (ES) to either LEO or GTO. 

The entire vehicle fleet starts at the Earth surface node and units are launched according to need, hence this node provides a supply of one unit of each tug, and the requisite number of crew vehicles (i.e. one \code{CSM} and one \code{LM} for each manned mission). In order to simulate the disposal of the \code{LM} after the surface crew has been transferred back to the \code{CSM}, the LLO node is modeled to have a demand for one \code{LM} every crew mission. Additionally, a fixed pre-calculated demand for \code{fLM} (LM fuel) is assigned to the LLO node for every crew mission as the amount required is fixed and is only consumed during the lunar descent/ascent operations. This information represents the demand/supply at nodes for use in Equation~\eqref{eq:balance}.

\subsection{Flow Transformation}
Commodities undergo transformation as they flow through the network, i.e., propellant is consumed in order to provide transportation across arcs. Since there are two kinds of propulsion available, the mass transformation matrix in Equation~\eqref{eq:xformation1} can take the two different forms listed below.
\setlist{nolistsep}
\begin{itemize}[noitemsep]
\item For high-thrust chemical propulsion, the $\Delta$v values gathered from the literature can be used with the Tsiolkovsky rocket equation to write the transformation constraint. The rocket equation can be written as:
\[ \displaystyle y^-_{ijve} = \exp \left(\frac{-\Delta v}{gI_\text{sp}} \right) y^+_{ijve}, \]
where $\Delta$v is the change in vehicle's velocity provided by its propulsion system, $I_\text{sp}$ is the corresponding specific impulse and $g$ is standard gravity. This above equation resembles the form of the mass transformation constraint derived in Equations~\eqref{eq:linear1}--\eqref{eq:xform_affine}, where $\displaystyle {}^1p = \exp \left(\frac{-\Delta v}{gI_\text{sp}} \right)$ and ${}^0p = 0$. 

\item For low-thrust solar electric propulsion, the trajectory analysis in Figure~\ref{fig:fit} reveals linear relationships between the final and initial masses on the arcs. Thus, using the values in Table \ref{tab:fit_mass}, the mass transformation constraint can again be implemented in the form of Equation~\eqref{eq:xform_affine}.
\end{itemize}

While the above discussion applies to transportation arcs within the event layers, it is assumed in the current case study that no commodity transformation occurs over holdover arcs as no propulsion is required in the case study to traverse these arcs. This implies that the matrix $\bm{B}_{ii}=\bm{I}_{k\times k}$, where $\bm{I}$ is an identity matrix and $k$ is the number of commodities.

Finally, it is important to note that in the case where the trajectory analysis yields a nonlinear relationship between the final and initial masses on the arcs, the nonlinear curve can be decomposed into piecewise linear segments and the MILP techniques described in Equation~\eqref{eq:pwl} be employed instead.

\subsection{Flow Concurrency and Bounds}
The concurrency constraints in Equations~\eqref{eq:concur1}--\eqref{eq:concur2} vary across the event layers as desired, but their general form for the case study campaign is listed below.
\setlist{nolistsep}
\begin{itemize}
\item Propellant being carried by a fixed-size spacecraft for its own transportation should not be exceed its maximum fuel capacity:
    \begin{equation}
    \begin{bmatrix} 1 &-C^{\rm fuel}_v  \end{bmatrix}_{ijv} \begin{bmatrix} \mathrm{propellant} \\ \mathrm{vehicle} \end{bmatrix}^+_{ijve} \leq 0, 
    \end{equation}
    \noindent where vehicle $v$ = \{\code{tug1}--\code{tug7}, \code{tug8}--\code{tug12}, \code{CSM}, \code{LM}\} and corresponding propellant = \{\code{fHIGH}, \code{fLOW}, \code{fCSM}, \code{fLM}\}. The fuel capacity of the corresponding vehicle is represented by $C^{\rm fuel}_v$. This form of the concurrency constraint applies to transportation arcs. 
    
    \item On the other hand, a different form of the constraint is applied when using a ``vehicle'' that is sized according to the propellant it is required to carry. In the case study, this only applies to the arcs that transport the crew to beyond TLI, which always use the launch vehicle's upper stage:
    \begin{equation}
    \begin{bmatrix} \hat{\varepsilon}_{\rm US} &  -1  \end{bmatrix}_{ijv} \begin{bmatrix} \mathrm{\code{fUS}} \\ \mathrm{\code{strUS}} \end{bmatrix}^+_{ijve} \leq 0, 
    \end{equation}
    \noindent where $\hat{\varepsilon}_{\rm US} = \varepsilon_{\rm US}/(1-\varepsilon_{\rm US})$. By accommodating \code{strUS} as a continuous variable (instead of fixing its size like the crew vehicles or the tugs), the savings in crew IMLEO achieved from a distributing the launch over smaller launch vehicles can be considered.
    
\item The sum of LM fuel and CSM fuel flowing across any arc (transportation or holdover) is not to exceed the sum of capacity of LM, CSM, and the droptank structural mass:
    \begin{equation}
    \begin{bmatrix} \hat{\varepsilon}_{\rm{LM}} & \hat{\varepsilon}_{\rm{CSM}} &  -1 &-\hat{\varepsilon}_{\rm{LM}} \cdot C_{\rm{LM}}^{\rm{fuel}} & -\hat{\varepsilon}_{\rm{CSM}} \cdot C_{\rm{CSM}}^{\rm{fuel}}  \end{bmatrix}_{ijv} \begin{bmatrix} \code{fLM} \\ \code{fCSM} \\ \code{strDtank} \\ \code{LM} \\ \code{CSM} \end{bmatrix}^+_{ijve} \leq 0, 
    \end{equation}
    \noindent where $ \hat{\varepsilon}_{\rm{LM}} = \varepsilon_{\rm{LM}}/(1-\varepsilon_{\rm{LM}})$, $ \hat{\varepsilon}_{\rm{CSM}} = \varepsilon_{\rm{CSM}}/(1-\varepsilon_{\rm{CSM}})$, with $\varepsilon_{\rm{CSM}}$ and  $\varepsilon_{\rm{LM}}$ representing the structural coefficients for CSM fuel and LM fuel respectively, and $C_{\rm{CSM}}^{\rm{fuel}}$ and  $C_{\rm{LM}}^{\rm{fuel}}$ are their fuel capacity. Note that the above equation assumes that LM and CSM tanks can be shared; if that assumption does not hold, each fuel tank type can be modeled as a different commodity type with separate concurrency constraints.
\end{itemize}

In addition to the above list, we also use the concurrency constraints together with the nonnegativity constraints in Equations~\eqref{eq:bound1}--\eqref{eq:bound2} to constrain certain variables to be zero when needed, i.e., $x^k_{ijve} = 0$ within the vector $\bm{x}^\pm_{ijve}$ when the commodity $k$ is not allowed on a given arc $(i,j,v,e)$. These constraints are set to reflect the following restrictions:
\setlist{nolistsep}
\begin{itemize}[noitemsep]
\item Cargo delivery can only be conducted by the available tugs. This means that crew vehicles (\code{LM}, \code{CSM} and \code{strUS}) are not allowed in the cargo delivery event layers. The opposite is also true where tugs (\code{tug1}--\code{tug12}) and tug propellants (\code{fHIGH} and \code{fLOW}) are not allowed in the crew layers.
\item The only commodities that can act as ``payloads'' in the cargo delivery layers are the droptank structural mass (\code{strDtank}), LM fuel (\code{fLM}) and CSM fuel (\code{fCSM}). This ensures that tugs are not piggybacking on each other.
\item In the cargo delivery events steps, each arc in a multi-graph set allows only one tug and its corresponding fuel type, along with desired payload commodities (i.e., \code{strDtank}).
\item The launch vehicle upper stage (and its corresponding fuel) is only allowed on the arcs that take the crew from LEO to TLI.
\end{itemize}

\subsection{Treating Time-Related Constraints}
The difficulty of including low-thrust transfer times encountered by existing space logistics optimization frameworks is dealt with by using Equation~\eqref{eq:tof_con} from the event-driven network optimization formulation. In the current case study, it is used to place bounds on the total time taken to set up the in-space propellant resupply chain (denoted henceforth by $T_{\rm cargo}$) and to constrain the total crew flight time across all missions ($T_{\rm crew}$) by splitting it into two conditions.

\subsubsection{Time of Flight}

The first step to deal with campaign-level time-constraints is to obtain the arc lengths of the transportation arcs enclosed within each event layer. The formulation for deriving the arc length from the relationships between the time of flight and the initial arc mass is shown in Equations~\eqref{eq:arcl}--\eqref{eq:pwl_t}. This section shows how they are specifically applied to the case study campaign.

If no vehicle transports any commodity across a particular pair of nodes in the event layer, then the time of flight over that arc is simply zero, i.e., $\mathrm{TOF}=0$. If a high-thrust vehicle is used to transport commodities across an arc in the case study campaign, the time of flight is fixed as:
\[ \mathrm{TOF} = {}^0q, \]
where ${}^0q$ takes a fixed value (specified in Table~\ref{tab:crewcosts} and \ref{tab:cptug}). On the other hand, the time of flight of the low-thrust arcs in the campaign is given by Equation~\eqref{eq:timefit}:
\[\mathrm{TOF}= {}^1q\cdot y^+_{ijve} + {}^0q, \]
where ${}^1q$ and ${}^0q$ are both constants that take values listed in Table~\ref{tab:fit_tof}. These above equations should be modified to reflect zero arc length when there is no flow across an arc. As explained previously, this is done by using $x^{\rm vehicle}_{ijve}$ which is the binary-valued component of the commodity flow vector $\bm{x}^+_{ijve}$ specifying the flow of vehicle $v$ over the arc ($i,j,v,e$). Thus the arc lengths of the transportation arcs in the current campaign are given by: 
\begin{equation}
\label{eq:arclength}    
\Delta t_{ijve}= 
\begin{cases}
{}^0q_{ijve}\cdot x^{\rm vehicle}_{ijve} \hquad \textrm{for\ high-thrust\ arcs}.\\
{}^1q_{ijve}\cdot y^+_{ijve} + {}^0q_{ijve}\cdot x^{\rm vehicle}_{ijve} \hquad \textrm{for\ low-thrust\ arcs}.
\end{cases}
\end{equation}

Conversely, if the arc length is dependent on the total mass of commodities flowing across the arc through a general nonlinear curve, it can instead be represented by decomposing into piecewise linear (PWL) segments using the procedure shown earlier in Equation~\eqref{eq:pwl_t}. 

\subsubsection{Cargo Delivery Duration}
The concept essential to incorporating the time of flight in the event-driven network model is the multi-graph, where an arc within the multi-graph with a particular index can only be traversed by the corresponding tug unit.
For example, arc \#1 in the multi-graph connecting two nodes $(i,j)$ within the event $e$ can only be used for transportation of commodities using \code{tug1}. Thus, if the set of transportation arcs that belong to event layer $e$ is $\mathbb{A}_e$, then the total time of flight $t_{ve}$ of tug $v$ over all the arcs within the layer $e$ is given by:
\begin{equation}
\label{eq:doc2}
t_{ve}  = \sum_{(i,j):(i,j,v)\in\mathbb{A}_e}\Delta t_{ijve} \quad \forall v \in \mathbb{V}_{\rm cargo} \quad \forall e \in \mathbb{E}_{\rm cargo}.
\end{equation}

\noindent where $\Delta t_{ijve}$ is the arc length derived earlier in Equation~\eqref{eq:arclength}. From here, the total duration $t_e$ of the event layer $e$ can be expressed as the maximum of sum of flight times of all tugs flowing in that layer:
\begin{equation}
\label{eq:doc3}
\max_{v\in \mathbb{V}_{\rm cargo}}\ (t_{ve}) = t_e\hquad.
\end{equation}

\noindent The above equality constraint in Equations~\eqref{eq:doc3} is converted into a linear inequality through the introduction of additional continuous variables $\widetilde t_{e}$ (one per event layer):
\begin{equation}
\label{eq:doc4}
 t_{ve} \leq \widetilde{t}_{e} \qquad \forall\ v \in \mathbb{V}_{\rm cargo}
\end{equation}

\noindent The desired bound on the duration of the campaign's cargo delivery phase can now be implemented as:
\begin{equation}
\label{eq:doc5}
\sum_{e \in \mathbb{E}_{\rm cargo} }\widetilde{t}_e \leq T_{\rm cargo}\hquad.
\end{equation}

Thus, by deriving the relationship between $\widetilde{t}_e$ and $\bm{x}^\pm_{ijve}$, this above Equation~\eqref{eq:doc5} can be written in a linear form and therefore can be integrated into the MILP formulation as Equation~\eqref{eq:tof_con}, where $\widetilde{t}_e$ corresponds to the auxiliary variable $\bm{\tau}$ in Equation~\eqref{eq:tof_con}. For the current case study campaign, this is done by condensing  Equation~\eqref{eq:arclength} into Equations~\eqref{eq:doc2}--\eqref{eq:doc5} as:
\begin{subequations}
\label{eq:tcargo}
\begin{align}
\sum_{(i,j):(i,j,v)\in\mathbb{A}_e}\bigg({{}^1q_{ijve}}\cdot y^+_{ijve} + {}^0q_{ijve}\cdot x^\text{vehicle}_{ijve}\bigg) - \widetilde{t}_e &\leq 0 \quad \forall\ v \in \mathbb{V}_{\rm cargo}  \quad \forall\ e \in \mathbb{E}_{\rm cargo}\hquad, \\
\sum_{e \in \mathbb{E}_{\rm cargo} }\widetilde{t}_e &\leq T_{\rm cargo}\hquad,
\end{align}
\end{subequations}

\noindent where ${}^1q_{ijve}=0$ for high-thrust arcs. (Note that $y^\pm_{ijve}$ can be written as a linear function of $\bm{x}^\pm_{ijve}$ using Equation~\eqref{eq:tot_mass}.) Equation~\eqref{eq:tcargo} is used to restrain the duration of the campaign's cargo delivery phase. 

\subsubsection{Crew Time of Flight}
Deriving a constraint on the total allowable crew time of flight across all missions is done in a similar manner and begins by simply calculating the total time of flight of the CSMs across all three missions. Every crew layer arc uses high-thrust propulsion, and thus corresponding arc lengths are given by:
\begin{equation}
\Delta t_{ijve} = {}^0q_{ijve}\cdot x^\text{vehicle}_{ijve}\hquad, \label{eq:toc1}
\end{equation}
\noindent where $v$ = \code{CSM} only and $e \in \mathbb{E}_{\rm crew}$. Due to the lack of alternative in-space crew vehicles in the analysis, the $\displaystyle \max_{v\in \mathbb{V}_{\rm cargo} }\ (t_{ve})$ relationship used in Equation~\eqref{eq:doc3} can be discarded and the time constraint for the crew phase is formulated directly as:
\begin{equation}
\sum_{e\in \mathbb{E}_{\rm crew}} \sum_{(i,j):(i,j,v)\in \mathbb{A}_e} \big({}^0q_{ijve}\cdot x^\text{vehicle}_{ijve} \big) \leq T_{\rm crew} \quad \forall\ v \in \mathbb{V}_{\rm crew}   \label{eq:tcrew}\hquad.
\end{equation}

\subsubsection{Summary of Campaign-Related Time Constraints}
Equations~\eqref{eq:tcargo} and \eqref{eq:tcrew} make up the vector of linear time-related constraints, which are the manifestations of Equation~\eqref{eq:tof_con} in the example event-driven network model. At this point, the derivation of the mathematical model is complete and can be supplied to a MILP optimizer. The decision variables available to the MILP-solver to optimize the event-driven network model for the case study campaign considered here are: $\left\lbrace \bm{x}^\pm_{ijve}, y^+_{ijve}, \bm{x}^\pm_{iie}, \widetilde{t}_e\right\rbrace $,
where $\bm{x}^\pm_{ijve}$ represents the commodity flow over each transportation arc, $y^+_{ijve}$ gives the total initial mass (i.e., outflow from the origin node) over each transportation arc, $\bm{x}^\pm_{iie}$ is the flow over each holdover arc, and $\widetilde{t}_e$ is the auxiliary variables introduced to denote the length of each cargo-related event layer. Additional $\lambda$ and $\lambda'$ variables can be added if PWL approximations are adopted for deriving the arc's trajectory performance metrics.

\section{Results \label{sec:results}}
This section presents the results of implementing the event-driven space logistics model for the design of in-space propellant resupply logistics for crew missions to the lunar surface. The MILP problem resulting from applying Equations~\eqref{eq:cost}--\eqref{eq:tot_mass} to the case study is implemented in MATLAB and solved using the Gurobi 6.5 solver on an Intel Core i7-4770 3.4 GHz platform. The results are compared against a baseline case consisting of no-refuel 3$\times$ Apollo missions by setting $T_{\rm cargo}$ to be zero and $T_{\rm crew} = 21$ days (7 days for each crew mission) within the mathematical model. This \textbf{baseline cost is IMLEO = 372.671 t} -- all optimal solutions obtained are expressed as percentage improvement over this value.

\subsection{Event-Driven Space Logistics Optimization Framework}
\begin{figure}[htb!]
\centering
\includegraphics[width=\textwidth]{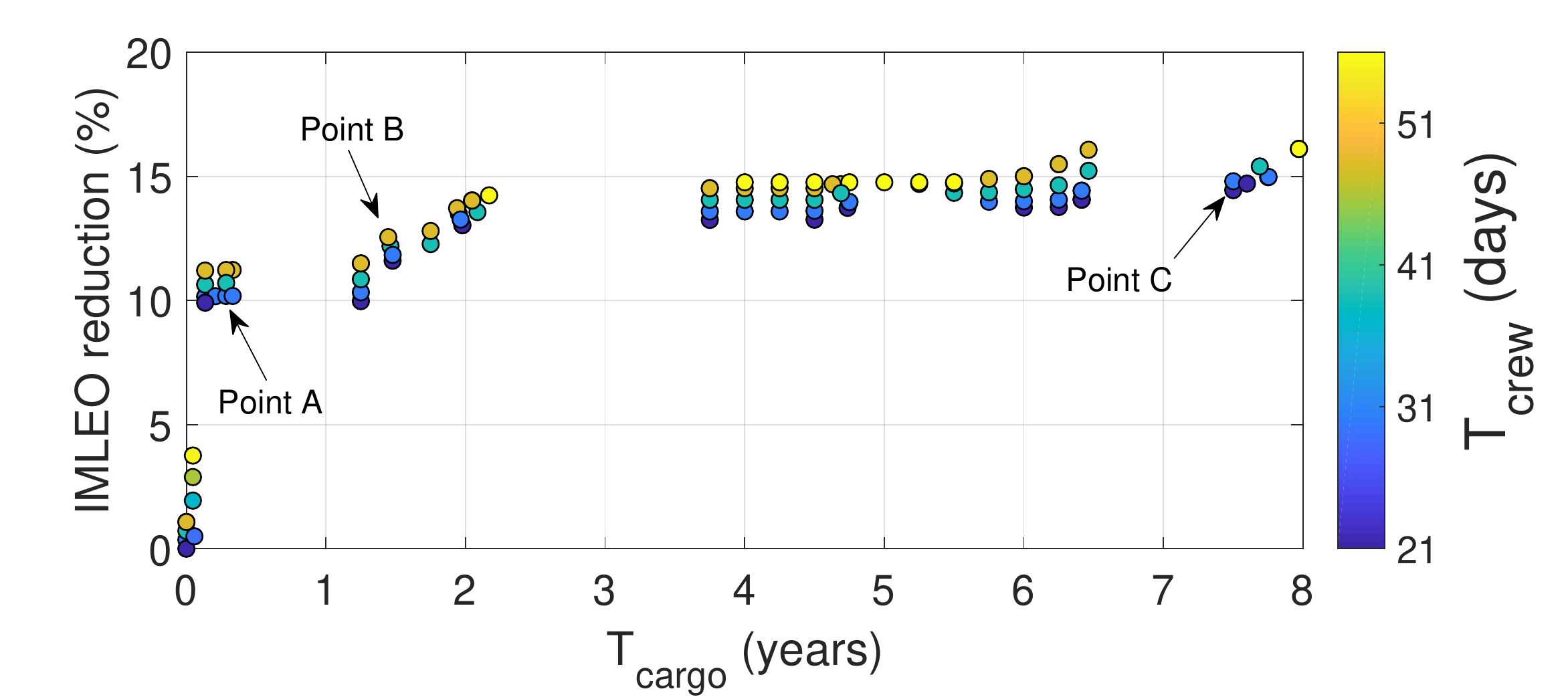}
\captionsetup{justification=centering}
\caption{Cost vs. time Pareto front for the case study obtained from applying the event-driven space logistics model. \label{fig:res1}}
\end{figure}

In order to demonstrate the true tradeoffs between cost and time, a Pareto front is manually generated for discrete values of maximum time bounds on crew flight times, and in steps of 4 months for maximum allowable time for cargo deliveries. This triple-objective front is displayed in Figure \ref{fig:res1} -- the color bar represents the total crew flight time across all three human missions. New options exploring in-space refueling architectures are available along this front. As expected, cost and time display a competing nature where more savings can be realized as the time constraint is relaxed. 

This Pareto front of optimal solutions is examined in detail by focusing on architectures depicting characteristics obtainable only by using the MILP-based event-driven space logistics optimization formulation developed in this work. These points are marked out in Figure~\ref{fig:res1}. The general trend is that architectures use chemical propulsion tugs only at the left-end of the Pareto front (Point A), and then proceed to use solar-electric propulsion tugs as the axis of time duration is broadened (Points B-C). Some salient features of the solution points on the Pareto front are:
\setlist{nolistsep}
\begin{enumerate}[noitemsep]
    \item \underline{Cargo relaying}: In some solutions (such as Points A-B), two tugs are used for cargo deliveries, but one of them is launched with minimal cargo to a way station in order to relay to a farther node the cargo that the other tug delivers. This behavior is further explored and explained in the following text.
    \item \underline{SEP tug reuse}: Each tug itself can be used up to a total of three times, by returning to its Earth parking orbit for refueling. The tug would, however, have to carry additional propellant to facilitate its return, apart from propellant required to make its cargo deliveries. However, such an architecture does not have to launch the dry mass of the tug into Earth orbit again, thus realizing some nontrivial savings. This solution is seen in Point C, where the architecture only launches the smallest SEP tug, but uses it twice. 
\end{enumerate}

The specific Pareto front obtained here is dependent on the modeling assumptions made, such as the vehicle fleet chosen for the problem and the $\Delta$v values. For example, the front may be smoother (i.e., without the wide gaps as seen) if more complex low-thrust trajectory models are made available. The front may also yield more solutions favoring SEP if more optimal low-thrust trajectory models are provided.

Let us now consider the solutions points marked in Figure \ref{fig:res1} in depth. 

\subsubsection{Point A}
\begin{figure}[htb!]
\centering
\begin{subfigure}{0.48\textwidth}
    \centering
    \includegraphics[clip,trim={0 0 12cm 0},height=0.37\textheight]{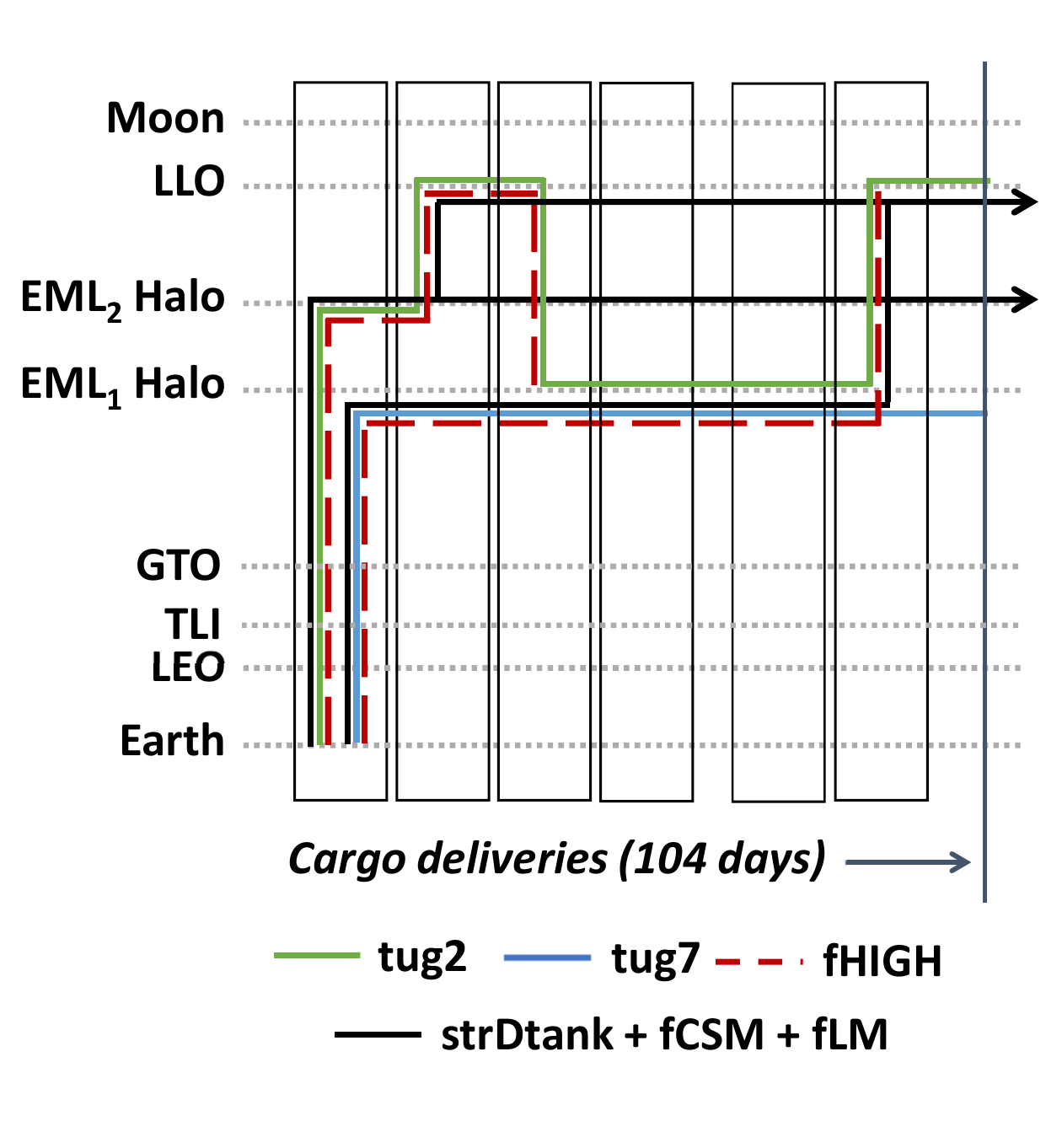}
    \caption{Cargo-related event layers (only E1 -- E6 shown) }
    \label{fig:A_cargo}
\end{subfigure}%
\begin{subfigure}{0.48\textwidth}
    \centering
    \includegraphics[clip,trim={15cm 0 0 0},height=0.37\textheight]{figures/PtA.pdf}
    \caption{Crew-related event layers (E13 -- E18) }
    \label{fig:A_crew}
\end{subfigure}
\caption{Commodity flows in solution point A from Figure~\ref{fig:res1}.}
    \label{fig:milpA}
\end{figure}

The commodity flow in solution Point A is depicted graphically in Figure~\ref{fig:milpA}, while its detailed breakdown is listed in Table~\ref{tab:pta}. The unique feature of this solution architecture is the relaying of cargo with different tugs in a very prominent manner.

\begin{table}[htb!]
\begin{adjustwidth}{-1cm}{-1cm}
\centering \footnotesize
\caption{Initial mass flow over each arc in solution point A from Figure \ref{fig:res1} \\ (integer variables specified in number of units used and continuous variables in kg).} \label{tab:pta}
\begin{tabular}{rrccccccccccc}
\hline\hline
Layer & Arc Name & \code{CSM} & \code{LM} & \code{strUS} & \code{fCSM} & \code{fLM} & \code{fUS} & \code{strDtank} & \code{fHIGH} & \code{tug2} & \code{tug7} & TOF \\
\hline \\[-8pt]
\multicolumn{13}{c}{\textit{Arcs within cargo delivery-related event layers (E1--E12)}} \\
1& ES  to LEO&0&0&0&14875&33140&0&4175&77780&1&1&0\\
1&LEO to L1&0&0&0&13735&29630&0&3771&68000&0&1&21\\
1 to 3&L1 to L1&0&0&0&13735&29630&0&3771&2988&0&1&-\\
1&LEO to L2&0&0&0&1139&3510&0&404&9780&1&0&17\\
1 to 2&L2 to L2&0&0&0&1139&3510&0&404&694&1&0&-\\
2&L2 to LLO&0&0&0&0&3510&0&305&694&1&0&27\\
2 to 3&LLO to LLO&0&0&0&0&3510&0&305&139&1&0&-\\
3&LLO to L1&0&0&0&0&0&0&0&139&1&0&28\\
3 to 10&L1 to L1&0&0&0&13735&29630&0&3771&2988&1&1&-\\
3 to 10&LLO to LLO&0&0&0&0&3510&0&305&0&0&0&-\\
10&L1 to LLO&0&0&0&13735&29630&0&3771&2988&1&0&28\\ 10 to 13 & L1 to L1 & 0 & 0 & 0 & 0 & 0 & 0 & 0 & 0 & 0 & 1 & -\\
10 to 13 & LLO to LLO & 0 & 0 & 0 & 13,735 & 33,140 & 0 & 4,076 & 0 & 1 & 0 & -\\
\\[-8pt]

\multicolumn{13}{c}{\textit{Arcs within crew flight-related event layers (E13--E18)}} \\13& ES  to LEO&1&1&4622&6711&0&35986&0&0&0&0&0\\
13&LEO to TLI&1&1&4622&6711&0&35986&0&0&0&0&0\\
13&TLI to LLO&1&1&0&6711&0&0&0&0&0&0&4\\
13 to 14&LLO to LLO&1&0&0&13735&22093&0&4076&0&0&0&-\\
14&LLO to L2&1&0&0&3364&0&0&0&0&0&0&3.5\\
14&L2 to  ES &1&0&0&1139&0&0&0&0&0&0&8.5\\
14 to 15&LLO to LLO&0&0&0&10372&22093&0&4076&0&0&0&-\\
15& ES  to LEO&1&1&4622&6711&0&35986&0&0&0&0&0\\
15&LEO to TLI&1&1&4622&6711&0&35986&0&0&0&0&0\\
15&TLI to LLO&1&1&0&6711&0&0&0&0&0&0&4\\
15 to 16&LLO to LLO&1&0&0&10372&11047&0&4076&0&0&0&-\\
16&LLO to  ES &1&0&0&5186&0&0&0&0&0&0&3\\
16 to 17&LLO to LLO&0&0&0&5186&11047&0&4076&0&0&0&-\\
17& ES  to LEO&1&1&4622&6711&0&35986&0&0&0&0&0\\
17&LEO to TLI&1&1&4622&6711&0&35986&0&0&0&0&0\\
17&TLI to LLO&1&1&0&6711&0&0&0&0&0&0&4\\
17 to 18&LLO to LLO&1&0&0&5186&0&0&4076&0&0&0&-\\
18&LLO to  ES &1&0&0&5186&0&0&0&0&0&0&3\\
\hline \\[-8pt]
\multicolumn{13}{c}{\normalsize Total IMLEO cost = 334.7268 t, T$_{\rm cargo}$ = 104 days, T$_{\rm crew}$ = 30 days} \\
\hline \hline
\end{tabular}
\end{adjustwidth}
\end{table}

Two chemical tugs are used to carry out all the cargo deliveries -- \code{tug2} (CP tug type 1), which is the smallest vehicle of its kind, and \code{tug7} (CP tug type 3), which is the largest high-thrust tug. \code{tug2} launches to the EML$_2$ halo with the propellant to resupply crews later on their return journey, drops off some of its payload, continues to LLO to again drop off a part of its payload, then returns to the EML$_1$ halo to wait for \code{tug7}. \code{tug7} then launches to LEO and travels to EML$_1$ with the rest of the LM fuel and CSM fuel required for crew resupply, along with required droptank structural mass, bundled as its cargo. Upon reaching the EML$_1$ node, this \code{tug7} transfers all of its cargo to the waiting \code{tug2}. More interestingly, \code{tug2} also receives additional tug propellant from the incoming \code{tug7}. This ultimately allows \code{tug2} to relay what was initially \code{tug7}'s cargo from EML$_1$ to LLO, where the droptank can then flow across holdover arcs to later refuel the incoming crew spacecraft (\code{LM} and \code{CSM}). In this solution architecture, the crew flight time is constrained to 30 days, and hence only one of the three crews visits the EML$_2$ halo for refueling on their return journey, while the other two crews are resupplied at LLO itself and return directly to Earth. 

Table~\ref{tab:pta} list the arcs clearly and tags them with the event layers that they belong to; as a reminder, layers 1 through 12 contain arcs that deliver cargo (i.e. the refuel droptanks), while the succeeding layers simulate crew journeys. The reader can also refer to Table \ref{tab:commodities} to understand the variable type used for each commodity; in general, all vehicles except the upper stage have discrete sizes and hence are modeled as binary variables. 

\FloatBarrier
\subsubsection{Point B}
The minimum length of the campaign's crew journey phase is 21 days, which derives from 7 days of flight time for each Apollo-style no-refuel crew mission. Longer crew flight times are permissible in other solutions, for example in those points marked by green or yellow color on the Pareto front. Consider Point B as an example of this case, where the crew flight is restricted to $\leq$50 days instead of just $\leq$30 days in solution point A. 
This relaxation in the crew flight time allows the crew to return via the EML$_2$ halo repeatedly, because the LLO$\rightarrow$ EML$_2\rightarrow$ Earth surface is cheaper than the direct return. This return through the EML$_2$ halo is then coupled with fuel resupply from predeployed droptanks to realize substantial savings ($\approx$12.55\%) over the baseline no-refuel case.

\begin{figure}[!htb]
\centering
\begin{subfigure}{0.48\textwidth}
    \centering
    \includegraphics[clip,trim={0 0 12cm 0},height=0.37\textheight]{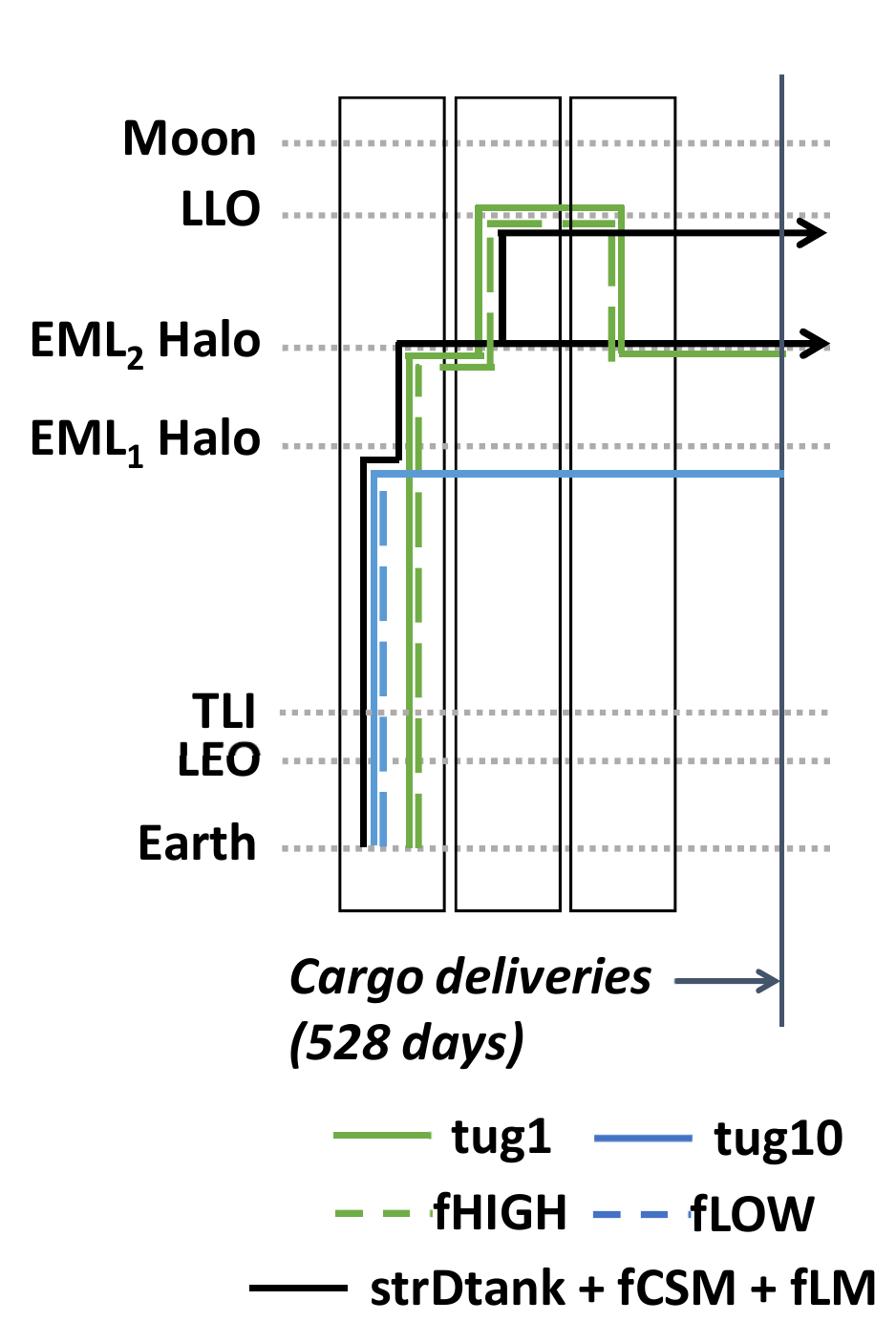}
    \caption{Cargo-related event layers (only E1--E3 shown) }
    \label{fig:B_cargo}
\end{subfigure}%
\begin{subfigure}{0.48\textwidth}
    \centering
    \includegraphics[clip,trim={11cm 0 0 0},height=0.37\textheight]{figures/PtB.pdf}
    \caption{Crew-related event layers (E13 -- E18) }
    \label{fig:B_crew}
\end{subfigure}
\caption{Commodity flows in solution point B from Figure~\ref{fig:res1}.}
    \label{fig:milpB}
\end{figure} 

Figure~\ref{fig:milpB} shows the campaign profile for the solution in Point B. This solution is showcased in some detail here because it exhibits the cooperative use of tugs for cargo relaying again, but with a CP tug and an SEP tug. Unlike Point A, the tug types used in this solution architecture differ in terms of propulsion systems, thus the tug propellant itself is not exchanged at in-space nodes. \code{tug1}  (smallest CP tug) launches to the EML$_1$ node without any cargo payload, but with a full tank of its own propellant. Here, it meets with \code{tug10} (mid-sized SEP tug) that brings the complete cargo supply, and then ferries these droptanks to LLO and further to EML$_2$. This trajectory route is possible due to the low-energy pathways that are exploited in this analysis.

\FloatBarrier
\subsubsection{Point C}
Points C lies farther on the right side of the Pareto front, where long campaign durations are allowed. However, its uniqueness lies in the reuse of a tug, which in this case is the smallest SEP tug available within the vehicle fleet. Since the crew flight time is constrained to 21 days in this solution, all crew propellant resupply has to occur at the LLO itself. Thus, \code{tug8} delivers this entire droptank cargo to LLO in two journeys. The time duration to set up this resupply chain, however, is very long because of the lower thrust that can be provided by \code{tug8}, as compared to \code{tug10} used in Point B. Due to the use of the smaller SEP tug for droptank delivery, this solution can achieve $\approx14.5\%$ savings over the no-refuel scenario.

\FloatBarrier

\subsection{Comparison with the State-of-the-Art Method}
The results obtained here are compared with those obtained from the state-of-the-art methodology available for campaign design of in-space refueling architecture \cite{jag_j1}. The approach in the referenced work uses a multiobjective genetic algorithm (MOGA) for campaign design, where different in-space architectures are represented by genetic sequences. Each gene pair in a fixed-length chromosome is used to indicate if the crew gets refueled at an intermediary node on their way to and back from the lunar surface (i.e. at EML$_1$, EML$_2$ or LLO), and the tug used to deliver the corresponding cargo. 

\begin{figure}[ht]
\centering
\includegraphics[width=\textwidth]{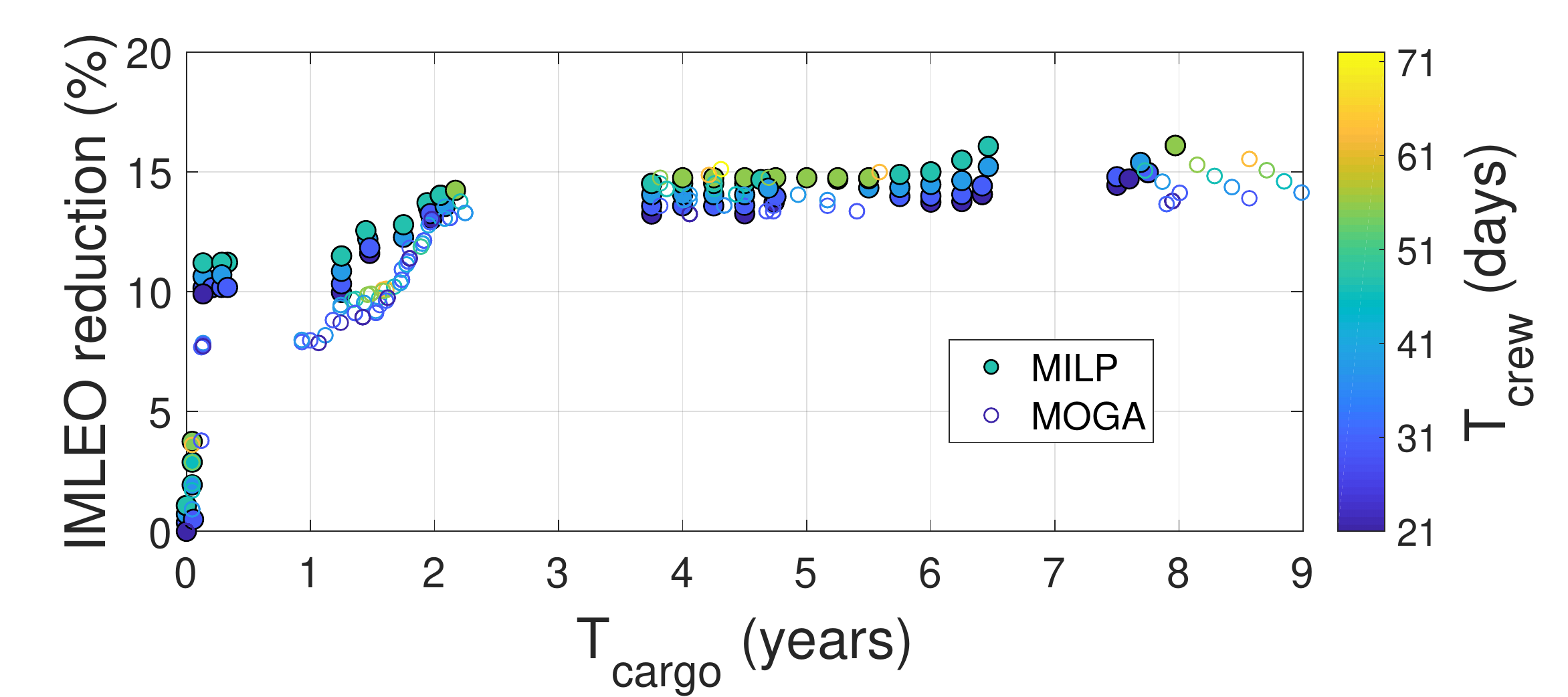}
\captionsetup{justification=centering}
\caption{Comparison of the Pareto front obtained from GA-based formulation and that obtained from the MILP-based method\label{fig:res2}}
\end{figure}

In order to showcase the advantages of using the new framework developed in the current work, the results are compared with those obtained from the existing MOGA-based formulation with the same problem settings (i.e. same fleet of vehicles and identical cost/performance models). The Pareto front obtained from the MOGA method is shown in Figure \ref{fig:res2} overlayed with the one obtained from the MILP-based method. The families of solutions that exist on this MOGA Pareto front are discussed in detail in the original work.

The highlight of this comparison is that the solutions obtained from the MILP-based event-driven space logistics optimization formulation realize higher overall cost savings than the earlier MOGA-based method. One of the core assumptions used in the MOGA-based formulation is that when the crew is refueled along its journey, it is given the exact amount of fuel to reach the next refuel droptank location or to complete the journey. Thus, the flexibility of carrying some of its own fuel is not afforded to the crew missions being resupplied. Furthermore, the MOGA-based formulation is deficient in the modeling of cargo/propellant relaying between tugs, and cannot exploit this powerful feature of the optimal solutions as generated by the event-driven space logistics model. Due to the lack of the rigorous optimization enabled by the MILP-based formulation, the quality of MOGA-based results is lower, i.e., the realizable IMLEO cost savings are decreased due to the limiting assumptions made on propellant refuel amount at way stations. By inherently eliminating restrictive assumptions, the proposed event-driven network approach is capable of exploring a substantially larger tradespace of solutions. 

As the time duration allowed for cargo deliveries is increased (i.e., moving towards the right along the Pareto front), more cargo can be delivered using SEP tugs. Due to this, at some point, almost all the propellant required to complete the crew journey is provided to it at in-space waystations, instead of carrying it along. In such cases, the SEP tugs are chosen to provide high-efficiency transportation, as long as the allowed duration is long enough. Hence, the solutions obtained from the event-driven space logistics model converge with those obtained from the MOGA-based formulation as we move further to the right along the time axis.  

Due to the more general problem formulation, the event-driven network method represents a substantial enhancement over the existing MOGA-based formulation. The increase in computational resources required to generate the Pareto front is not significant. While a single solution with the MILP-based event-driven network method takes only $\sim$3 seconds on the platform being used, the manual generation of the entire Pareto front takes about 15 minutes. The MOGA, on the other hand, completes its run in little over a minute (200 generations using a population size of 300 candidates). It must also be noted that the MOGA approach is stochastic in nature and its solutions are often sensitive to the initial population of the candidate solutions used. This is in contrast with the deterministic nature and guaranteed optimality (with respect to the given surrogate models) of the results obtained from the MILP formulation, thus highlighting another advantage of the MILP-based event-driven space logistics model.

\section{Conclusion}
This paper proposes a new event-driven space logistics model to optimize a multi-mission campaign with both high-thrust and low-thrust propulsion options. The challenge in this problem is that the transfer performance metrics (in terms of final-to-initial mass ratio and time of flight) for low-thrust arcs can vary with the mass of their commodity flows. In order to incorporate this flow-dependent arc length, the proposed event-driven space logistics optimization framework considers the time dimension through the use of event-based time steps, instead of discrete-sized time-steps as has been done previously. In addition, the formulation for the event-driven space logistics optimization framework considers the trajectory performance models using general surrogate models (affine or nonlinear) and incorporating them into a mixed-integer linear programming (MILP) formulation with reasonable approximation so that the space campaign design problem can be solved computationally efficiently. 

The utility of the developed framework is illustrated by applying it to the design of a campaign to establish a cislunar crew propellant resupply chain using tugs with high-thrust propulsion and low-thrust propulsion cooperatively -- this propellant resupply chain serves three Apollo-type crew missions. Low-energy pathways available in the Earth-moon system are exploited for cislunar cargo transport. With the model formulated for the case study campaign, optimal solutions (in terms of initial-mass-in-LEO cost) found for a range of time constraints for cargo and crew are used to populate the Pareto front. Leveraging in-space refueling options, the benefit of harnessing the full potential of decoupling the flight-time requirement for crew from that for cargo deliveries is explored. The most notable feature of the solutions on the Pareto front is the optimal relaying of cargo between different sets of tugs to collectively set up the in-space crew propellant resupply chain. Another highlight is the tug reuse when long cargo delivery times are allowed; tugs are considered reused when they return to their Earth parking orbit to be refueled and collect more cargo. 

In addition, the event-driven space logistics solutions for the case study are also compared with architectures suggested by the state-of-the-art method for solving the same case study as considered here. The current model eliminates certain limiting assumptions from the existing model, thus making the problem formulation more general and exploring a larger tradespace of campaign architectures. The insignificant increase in the computational effort to populate the Pareto front is worthwhile because of the higher quality of solutions that can be obtained.

The campaign profiles derived from the event-driven space logistics model not only reduce costs but can also be used for recommending roadmaps for technology development, such as in-space propellant storage and transfer, deep-space rendezvous, and solar electric propulsion tugs for cargo
delivery. Through the use of the MILP-based event-driven network optimization framework presented in this paper, architectural decisions to deploy these technologies can be supported through rigorous mathematical modeling. 

As a final note, the event-driven space logistics model is developed with the motivation of enabling embedded propulsion technology trades for optimizing space logistics problems. This developed optimization method can trade off cost, time and technology in an automated manner, and thus can be effectively applied in other complex logistics situations involving flow-dependent arc lengths, where the time distribution between missions segments is not known beforehand.

\section*{Acknowledgment}
We appreciate Prof. Natashia Boland's feedback on this work.

\section*{Appendix A: Transfer Performance Metrics for High-Thrust Arcs}
\setcounter{table}{0}
\renewcommand{\thetable}{A\arabic{table}}

\begin{table}[htb!]
\centering
\caption{Crew trajectory performance metrics \label{tab:crewcosts}}
\begin{threeparttable}
\begin{tabular}{rccl}
\hline \hline
Arc & $\Delta$v, km/s & TOF, days & Impulse provided by \\ \hline \\[-8pt]
\multicolumn{4}{c}{\textit{Forward direction arcs\tnote{**}}} \\
LEO to TLI\tnote{*} & 3.306 & - & Launch vehicle upper stage \\
TLI to LLO\tnote{*} & 0.976 & 4 & CSM \\[2pt]
LEO to TL1I\tnote{*} & 3.200 & - & Launch vehicle upper stage \\
TL1I to L1 \cite{man2eml3} & 0.570 & 5 & CSM \\
L1 to LLO \cite{farq}  & 0.750 & 3 & CSM \\[2pt]
LEO to TL2I\tnote{*} & 3.400 & - & Launch vehicle upper stage \\
TL2I to L2 \cite{man2eml3} & 0.275 & 8.5 & CSM \\
L2 to LLO \cite{farq} & 0.750 & 3.5 & CSM \\[3pt]
\multicolumn{4}{c}{\textit{Return direction arcs\tnote{**}}} \\
LLO to ES\tnote{*} & 1.091 & 3 & CSM \\[2pt]
LLO to L1 \cite{farq} & 0.750 & 3 & CSM \\
L1 to ES \cite{man2eml3} & 0.499 & 10 & CSM \\[2pt]
LLO to L2 \cite{farq} & 0.750 & 3.5 & CSM \\
L2 to ES \cite{man2eml3} & 0.275 & 8.5 & CSM \\
\hline \hline
\end{tabular}
  \begin{tablenotes}
    \item TLI = Translunar injection, TL1I = Translunar injection to go to L$_1$ halo, TL2I = Translunar injection to go to L$_2$ halo.
    \item[*] Based on/ calculated from Apollo 17 values \cite{apollo}.
    \item[**] Lunar descent $\Delta$v = 2.123 km/s and lunar ascent $\Delta$v = 2.239 km/s.
    \item
  \end{tablenotes}
\end{threeparttable}
\end{table}

\begin{table}[htb!]
\centering
\caption{High-thrust CP tug trajectory performance metrics \label{tab:cptug}}
\begin{tabular}{rcc}
\hline \hline
Arc & $\Delta$v, km/s & TOF, days \\ \hline 
LEO to/from L1 \cite{mingtao} & 3.375 & 21 \\
L1 to/from LLO \cite{jag_j2} & 0.259 & 28 \\[2pt]
LEO to/from L2 \cite{mingtao} & 3.336 & 17 \\
L2 to/from LLO \cite{jag_j2} & 0.375 & 27 \\
\hline \hline
\end{tabular}
\end{table}

\section*{Appendix B: Vehicle Fleet}
\setcounter{table}{0}
\renewcommand{\thetable}{B\arabic{table}}
\FloatBarrier

\begin{table}[htp]
\centering
\caption{Specifications of crew vehicles used in analysis \label{tab:crewveh}}
\begin{threeparttable}
\begin{tabular}{lccc}
\hline \hline
    & Upper Stage (US) & CSM & LM \\ \hline
 Dry mass, t & Sized according &12.2 & 5.8 \\
    & to TLI & & \\
 Propellant capacity, t & - & 31 & 12 \\
 $I_\text{sp}$, s & 421 & 314 & 311 \\
 Structural coefficient of fuel ($\varepsilon$) & 0.1138\tnote{+} & 0.08\tnote{++} & 0.08\tnote{++}\\
\hline \hline
\end{tabular}
  \begin{tablenotes}
    \item[+] Calculated based on $\sim$14t dry mass of Saturn V upper stage and its 123t propellant capacity \cite{apollo}.
    \item[++] Used to size the required droptanks.
  \end{tablenotes}
\end{threeparttable}
\end{table}
\FloatBarrier

\section*{Appendix C: Linear Trends in Low-Thrust Transfer Performance Metrics}
Consider a simple case of very low continuous thrust between two circular coplanar orbits. Using tangential thrust in this case, the time of flight can be shown to be inversely proportional to the spacecraft's thrust-acceleration \cite{orbitalbook}:
\begin{equation*}
t_f = \frac{1}{\Gamma} \left(\sqrt{\frac{\mu}{a_i}}-\sqrt{\frac{\mu}{a_f}}\right)
\end{equation*}
\noindent where $\Gamma$ is the thrust acceleration, $a_i$ and $a_f$ are the semi-major axes of the initial and final orbits respectively, and $t_f$ is the time of flight. 

This relationship is similar to the almost-linear trends visible in Figure \ref{fig:fit}, although the transfer does not involve circular orbits. In this case, low-thrust low-energy transfers are designed between two points of fixed energies (because manifold target points are fixed) through the use of Q-law. Thus as payload being carried increases, the thrust-acceleration decreases, which then drives up the time of flight in an almost linear manner as well.    

\bibliography{references}

\end{document}